\documentclass[journal]{IEEEtran}
\usepackage{amssymb}
\usepackage{amsthm}
\usepackage[cmex10]{amsmath}
\usepackage{graphicx}
\usepackage{colortbl}
\usepackage{cite}

\newtheorem{theorem}{Theorem}

\newtheorem{remark}{Remark}
\newtheorem{lemma}{Lemma}
\newtheorem{assumption}{Assumption}
\newtheorem{definition}{Definition}

\newtheorem{observation}{Observation}

\ifCLASSINFOpdf
\else
\fi
\begin{document}
%
\title{On Irregular Linear Quadratic Control: Stochastic Case}
%
%
%

\author{Huanshui~Zhang~ and ~Juanjuan~Xu 
\thanks{This work is supported by the National Natural
Science Foundation of China under Grants 61633014, 61573221.  }
\thanks{H. Zhang is with School of Control
Science and Engineering, Shandong University, Jinan, Shandong, P.R.China 250061.
        {\tt\small hszhang@sdu.edu.cn}}
\thanks{J. Xu is with School of Control Science and Engineering, Shandong University, Jinan, Shandong, P.R. China 250061.
        {\tt\small jnxujuanjuan@163.com}}
}

\maketitle

\begin{abstract}
As it is popular known, Riccati equation is the key basic tool for optimal control in the modern control theory. The solvability conditions of optimal control, stabilization conditions and controller design are all based on the Riccati equation. However, these results highly depends on a strictly assumption that the Riccati equation is regular. If the Riccati equation is irregular, the controller could not be derived from the equilibrium condition. This paper is concerned with the general stochastic linear quadratic (LQ) with irregular Riccati equation.  Different from the classical control theory for regular LQ problems, a new approach of `multi-layer optimization' is proposed. With the approach, we show that different controller entries of irregular-LQ controller need to be derived from different  equilibrium  conditions and a specified terminal constraint condition in different layers, which is much different from the classical regular LQ control where all the controller entries can be obtained from equilibrium  condition in one layer based on regular Riccati equation. The presented results clarify the differences of open-loop control from closed-loop control in the aspects of solvability and controller design and also explores in essentially the differences of regular control from irregular control. Several examples are presented to show the effectiveness of the proposed approach.
\end{abstract}

\begin{IEEEkeywords}
Indefinite LQ control; Controllability; Riccati equation; Open-loop solvability; Closed-loop solvability.
\end{IEEEkeywords}

%
\IEEEpeerreviewmaketitle

\section{Introduction}

Consider linear control systems governed by  It\^{o} stochastic differential equation as
\begin{eqnarray}
dx(t)&=&[A(t)x(t)+B(t)u(t)]dt+[\bar{A}(t)x(t)\nonumber\\
&&+\bar{B}(t)u(t)]dw(t),~x(t_0)=x_0,\label{s1}
\end{eqnarray}
where $x\in R^n$ is the state, $u\in R^m$ is the control input. $w(t)$ is a standard one-dimension Brownian motion.
The matrices $A(t),B(t),\bar{A}(t),\bar{B}(t)$ are deterministic matrices with appropriate dimension.

The cost function is given by
\begin{eqnarray}
J(t_0,x_0;u)&=&E\int_{t_0}^T[x'(t)Q(t)x(t)+u'(t)R(t)u(t)]dt\nonumber\\
&&+Ex'(T)Hx(T),\label{s2}
\end{eqnarray}
where $Q(t),R(t),H$ are symmetric matrices with appropriate dimensions. Denote
\begin{eqnarray}
\mathcal{U}[t_0,T]&=&\{u(t),t\in[t_0,T]|u(t)~\mbox{is}~\mathcal{F}_{t}~\mbox{adapted},\nonumber\\
&&~~
E\int_0^T\|u(t)\|^2dt<\infty\}.\nonumber
\end{eqnarray}
Then the optimal linear quadratic (LQ) control problem is stated as following:\\

\textbf{Problem (Indefinite LQ):} For any given initial pair $(t_0,x_0),$ find a $u^*\in \mathcal{U}[t_0,T]$ such that
\begin{eqnarray}
J^*(t_0,x_0;u^*)=\min_{u(\cdot)\in\mathcal{U}[t_0,T]}J(t_0,x_0;u)\doteq V(t_0,x_0).\nonumber
\end{eqnarray}
The following assumption of {\em convexity} is necessary for the indefinite LQ control.
\begin{assumption}\label{a1}
\begin{eqnarray}
J(t_0,0;u)\geq0.\nonumber
\end{eqnarray}
\end{assumption}

The research of the above Problem  has been thoroughly studied since 1950's of last century \cite{Wonham},\cite{bismut}, \cite{bellman}, \cite{zhou}. A brief summary for the related works is given below.
\begin{itemize}
\item In the case of $\bar{A}(t)=0$ and $\bar{B}(t)=0$, (1) is a deterministic system, Problem (Indefinite LQ) is reduced to the deterministic optimal LQ control problem. The deterministic LQ can be traced back to 1950's of last century, please see \cite{bellman1} in 1958, \cite{kalman} in 1960, \cite{letov} in 1961, \cite{anderson} in 1990 and references therein.  When the LQ is regular (i.e., the associated Riccati equation is regular), the problem can be solved elegantly via the Riccati equation; see \cite{anderson} for a thorough study of the Riccati equation approach.



\item Stochastic LQ problems, i.e., $\bar{A}(t)\not=0$ and/or $\bar{B}(t)\not=0$ in (1), pioneered by \cite{Wonham}, has been widely studied, see \cite{29}-\cite{bismut} and references therein. In the literature it is typically assumed that the cost function has a positive definite weighting matrix for the control term,
and a positive semi-definite weighting matrix for the state term \cite{davis}. A  new type Riccati equation (generalized Riccati equation) has been applied to characterize the
solvability of the stochastic LQ control problem \cite{Wonham2},\cite{bismut} under the positive/seimi-positive assumptions. Another significant progress on stochastic control was made by \cite{zhou} in 1998. It was shown that a stochastic LQ problem with indefinite $Q$ and $R$ may still be well-posed. The controller could be designed in terms of state feedback via a stochastic/general Riccati equation. Follow-up research on indefinite stochastic LQ control and the general Riccati equation have been carried out in \cite{zhou2} and references therein.

\end{itemize}
However, it is noted that all the aforementioned results on optimal LQ control are based on the regularity of a Riccati equation. So a question is naturally arised: What is the case when the Riccati equation is irregular. To answer the question, we first recall the Maximum Principle for the above LQ control problem.

\begin{lemma}\label{lemMP}
(Maximum Principle) Under the convexity Assumption \ref{a1}, there exists a solution to Problem (Indefinite LQ) if and only if there exists $u(t)$ such that
\begin{eqnarray}
0&=&R(t)u(t)+B'(t)p(t)+\bar{B}'(t)q(t),\label{3}
\end{eqnarray}
where $p(t),q(t)$ satisfy the following forward backward stochastic differential equations (FBSDEs):
  \begin{eqnarray}
dx(t)&=&[A(t)x(t)+B(t)u(t)]dt+[\bar{A}(t)x(t)\nonumber\\
&&+\bar{B}(t)u(t)]dw(t),\label{1}\\
dp(t)&=&-[A'(t)p(t)+\bar{A}'(t)q(t)\nonumber\\
&&+Q(t)x(t)]dt+q(t)dw(t), p(T)=Hx(T). \label{2}
\end{eqnarray}
\end{lemma}
In usually, it is assumed that the relationship of costate $p(t)$ and state $x(t)$ is homogeneous, i.e.,
\begin{eqnarray}
p(t)= P(t)x(t),\label{Intr1}
\end{eqnarray}
where $P(t)$ is the solution to the following Riccati equation
\begin{eqnarray}
0&=&\dot{P}(t)+A'(t)P(t)+\bar{A}'(t)P(t)\bar{A}(t)+P(t)A(t)\nonumber\\
&&+Q(t)-\Gamma_0'(t)\Upsilon_0^{\dag}(t)\Gamma_0(t),~~P(T)=H, \label{r1}
\end{eqnarray}
while
\begin{eqnarray}
\Upsilon_0(t)&=&R(t)+\bar{B}'(t)P(t)\bar{B}(t),\label{Intr3}\\
\Gamma_0(t)&=&B'(t)P(t)+\bar{B}'(t)P(t)\bar{A}(t).\label{Intr4}
\end{eqnarray}Then the equilibrium condition (\ref{3}) becomes
\begin{eqnarray}
0=\Upsilon_0(t)u(t) +\Gamma_0(t)x(t). \label{Intr2}
\end{eqnarray}
If the Riccati equation is irregular, i.e., $Range \Big(\Gamma_0(t)\Big)\not\subseteq Range \Big(\Upsilon_0(t)\Big)$, it is clear that the controller $u(t)$ can not be solved from the above equilibrium (arbitrary $x(t)$). However, it does not implies that the LQ control is unsolvable. In fact, the following example shows that there still exists optimal controller even if the Riccati equation is irregular.

Example \cite{sunjingrui}: Consider the system
\begin{eqnarray}
dx(t)&=&[u_1(t)+u_2(t)]dt+[u_1(t)-u_2(t)]dw(t),\nonumber\\
&&~~~~~~~~~~~t\in[t_0,1]\label{r22}
\end{eqnarray}
and the cost function
\begin{eqnarray}
J(0,u)=Ex(1)^2.\nonumber
\end{eqnarray}
The optimal cost is obviously $0.$ The corresponding open-loop optimal controller is
\begin{eqnarray}
u(t)
&=&-\left[
      \begin{array}{c}
        1 \\
        1 \\
      \end{array}
    \right]\frac{x_0}{2(1-t_0)},\nonumber
\end{eqnarray}
and the closed-loop optimal controller is given by
\begin{eqnarray}
u(t)
&=&\left[
      \begin{array}{c}
        1 \\
        1 \\
      \end{array}
    \right]\frac{x(t)}{2(t-1)}.\nonumber
\end{eqnarray}
On the other hand, the solution to (\ref{r1}) is $P(t)=1.$ Thus $\Gamma_0(t)=\left[
                                                                            \begin{array}{c}
                                                                              1 \\
                                                                              1 \\
                                                                            \end{array}
                                                                          \right]
$ and $\Upsilon_0(t)=\left[
                     \begin{array}{cc}
                       1 & -1 \\
                       -1 & 1 \\
                     \end{array}
                   \right].
$ This implies that the Riccati equation is irregular. However, we noted that the problem is solvable.

It should be pointed out that the irregular LQ problem has been studied in earlier works \cite{zhou2}, \cite{sunjingrui}.  \cite{sunjingrui} pioneered the study on irregular LQ problem, focusing on the closed-loop and open-loop solvability for the stochastic control. \cite{sunjingrui} showed  that the solvability obtained in previous work \cite{zhou2} is only applicable to the case of closed-loop control, and open-loop solvability is different from closed-loop solvability. In particular when the Riccati equation is irregular, by constructing a minimizing sequence, the open-loop solvability is equivalent to the weak convergence and strong convergence of the minimizing sequence respectively.

%
%

In this paper, following-up the previous works \cite{zhou2,sunjingrui}, we will further study the optimal LQ control problem to clarify the above questions. Firstly, it is shown that the irregularity of Riccati equation implies the
nonhomogeneous relationship of the state and costate, but not the non-solvability of the LQ control. Both the open-loop and closed-loop solvability conditions are depicted by the Range conditions and terminal constrain condition (i.e., $Mx(T)=0$ with specified matrix $M$). Secondly, different from the classical regular-LQ where the controller are obtained by using the equilibrium  condition and the Riccati equation in one layer, the Irregular-LQ problems need more layers in order to obtain the controller entries. In other words, different controller entries need to be derived with different equilibrium conditions and terminal condition. Thirdly, the difference of open-loop and closed-loop control lies in that the former is to seek the entries of controller in open-loop form to
satisfy the terminal constraint condition while the latter to seek entries of controller in closed-loop form to
satisfy  the terminal constraint condition.

%
%
%

The remainder of the paper is organized as follows. Section II presents the problem and illustrates some preliminaries used
in the derivation of the controller. The regular and irregular solutions are given in Section III and IV respectively.
Section V shows the open-loop and closed-loop solutions. An example is illustrated in Section VI.
Some concluding remarks are given in Section VII.

The following notations will be used throughout this paper: $R^n$
denotes the family of $n$ dimensional vectors. $x'$ means the
transpose of $x.$ A symmetric matrix $M>0\ (\geq 0)$ means
strictly positive definite (positive semi-definite). $(\Omega, \mathcal{F}, \mathcal{P}, \mathcal{F}_t|_{t\geq 0})$ is a complete
stochastic basis so that $\mathcal{F}_0$ contains all P-null
elements of $\mathcal{F},$ and the filtration is generated by the
standard Brownian motion $\{w(t)\}_{t\geq 0}.$ We also introduce the following set:
\begin{eqnarray}
L_{\mathcal{F}}^2(0,T;R^m)\hspace{-2mm}&=&\hspace{-2mm}\{\varphi(t)_{t\in[t_0,T]}
\mbox{~is~an~}\mathcal{F}_t~\mbox{adapted~stochastic}\nonumber\\
&&\mbox{~process}~s.t.~
E\int_0^T\|\varphi(t)\|^2dt<\infty\}.\nonumber
\end{eqnarray}

\section{Definitions and Preliminaries }

%
%
We firstly introduce the following definitions.
\begin{definition}
Problem (Indefinite LQ) is said to be finite at initial pair $(t_0,x_0)\in[0,T]\times R^n,$ if
\begin{eqnarray}
V(t_0,x_0)>-\infty.\label{def1}
\end{eqnarray}
Problem (Indefinite LQ) is said to be finite at $t_0\in[0,T]$ if (\ref{def1}) holds for all $x_0\in R^n.$
Problem (Indefinite LQ) is said to be finite if (\ref{def1}) holds for all $(t_0,x_0)\in[0,T]\times R^n.$
\end{definition}

\begin{definition}
An element $u^*\in \mathcal{U}[t_0,T]$ is called an open-loop optimal control of Problem (Indefinite LQ) for the initial pair $(t_0,x_0)\in[0,T]\times R^n,$ if
\begin{eqnarray}
J^*(t_0,x_0;u^*)\leq J(t_0,x_0;u), \forall u(\cdot)\in\mathcal{U}[t_0,T].\label{def2}
\end{eqnarray}
If an open-loop optimal control (uniquely) exists for $(t_0,x_0)\in[0,T]\times R^n,$ Problem (Indefinite LQ) is said to be (uniquely) open-loop solvable at $(t_0,x_0)\in[0,T]\times R^n.$ Problem (Indefinite LQ) is said to be (uniquely) open-loop solvable at $t_0\in[0,T)$ if for the given $t_0,$ (\ref{def2}) holds for all $x\in R^n.$ Problem (Indefinite LQ) is said to be (uniquely) open-loop solvable if it is open-loop solvable at all $(t_0,x_0)\in[0,T]\times R^n.$
\end{definition}

\begin{definition} \cite{Oksendal}
An element $u^*(\cdot)=K(\cdot)x^*(\cdot)\in\mathcal{U}[t_0,T]$  is called a closed-loop optimal control of Problem (Indefintie LQ) on $[t_0,T]$ if
\begin{eqnarray}
J^*(t_0,x_0;K(\cdot)x^*(\cdot))\leq J(t_0,x_0;u), ~\forall x_0\in R^n,\nonumber\\~\forall u(\cdot)\in\mathcal{U}[t_0,T].\label{d3}
\end{eqnarray}
where $x^*(\cdot)$ is the solution to the following closed-loop system:
\begin{eqnarray}
dx^*(t)&=&[A(t)+B(t)K(t)]x^*(t)dt+[\bar{A}(t)\nonumber\\
&&+\bar{B}(t)K(t)]x^*(t)dw(t),~x(t_0)=x_0.\label{d4}
\end{eqnarray}
If a closed-loop optimal controller exists on $[t_0,T],$ Problem (Indefinite LQ) is said to be (uniquely) closed-loop solvable
on $[t_0,T].$ Problem (Indefinite LQ) is said to be (uniquely) closed-loop solvable if it is (uniquely) closed-loop solvable
on any $[t_0,T].$
\end{definition}

\begin{remark}
The above definition of closed-loop solvability is from \cite{Oksendal}, which is different from  the one given in \cite{sunjingrui}.
In fact, in \cite{sunjingrui}, it is assumed that $K(\cdot)\in L^2([0,T];R^{m\times n}).$  While in this paper,
it is only assumed that $K(\cdot)x^*(\cdot)\in \mathcal{U}[t_0,T]$ which allows the
case that $K(\cdot)\not\in L^2([0,T];R^{m\times n})$, please see \cite{Oksendal} for details.
\end{remark}

We now recall the Moore-Penrose inverse of a matrix which will be  used this paper repeatedly. From \cite{pinv}, for
a given matrix $M\in R^{n\times m}$, there exists a unique matrix in
$R^{m\times n}$ denoted by $M^{\dag}$ such that
\begin{eqnarray}
MM^{\dag}M=M,~M^{\dag}MM^{\dag}=M^{\dag},\nonumber\\
(MM^{\dag})¡¯=MM^{\dag},~(M^{\dag}M)'=M^{\dag}M.\nonumber
\end{eqnarray}
The matrix $M^{\dag}$ is called the Moore-Penrose inverse of $M$. The following lemma is from \cite{rami}.
\begin{lemma}\label{lem}
Let matrices $L, M$ and $N$ be given with appropriate size. Then, $LXM=N$ has a solution $X$ if and only if $LL^{\dag}NMM^{\dag}=N$. Moreover, the solution of $LXM=N$ can be expressed as
$X=L^{\dag}NM^{\dag}+Y-L^{\dag}LY MM^{\dag},$ where $Y$ is a matrix with appropriate size.
\end{lemma}
In particular, let $M=I$, we have $LX=N$ has a solution if and only
if $LL^{\dag}N = N$. This is also equivalent to $Range(N)\subseteq
Range(L)$ where $Range(N)$ is the range of $N.$

\section{Optimal control with regular Riccati equation}

For the completeness of presentation and comparisons with classical results, we will first consider the classical regular-LQ problem. Recalling the Maximum Principle in Lemma \ref{lemMP}, we note that the key to solve Problem (Indefinite LQ) is to obtain the solution to (\ref{3})-(\ref{2}).
\\

Different from the previous works as in (\ref{Intr1}), we assume without loss of generality that
\begin{eqnarray}
\Theta(t)\doteq p(t)-P(t)x(t),\label{21}
\end{eqnarray}
where $P(t)$ obeys Riccati equation (\ref{r1}). It is obvious that $P(T)=H$ and $\Theta(T)=p(T)-P(T)x(T)=0$. We also assume without loss of generality that
\begin{eqnarray}
d\Theta(t)=\hat{\Theta}(t)dt+\bar{\Theta}(t)dw(t), \label{zj1}
\end{eqnarray}
where $\hat{\Theta}(t)$ and $\bar{\Theta}(t)$ are to be determined.  Applying It\^{o}'s formula to (\ref{21}) yields
\begin{eqnarray}
dp(t)&=&\dot{P}(t)x(t)dt+P(t)\big[A(t)x(t)
+B(t)u(t)\big]dt\nonumber\\
&&+P(t)\big[\bar{A}(t)x(t)+\bar{B}(t)u(t)\big]dw(t)\nonumber\\
&&+\hat{\Theta}(t) dt+\bar{\Theta}(t)dw(t).\label{zz4}
\end{eqnarray}
Using (\ref{21}), we rewrite (\ref{2}) as
\begin{eqnarray}
dp(t)&=&-\big[A'(t)P(t)x(t)+A'(t)\Theta(t)
+\bar{A}'(t)q(t)\nonumber\\
&&+Q(t)x(t)\big]dt+q(t)dw(t).\label{11}
\end{eqnarray}
With a comparison of (\ref{zz4}) and (\ref{11}), it follows that
\begin{eqnarray}
q(t)&=& P(t)\big[\bar{A}(t)x(t)+\bar{B}(t)u(t)\big]+\bar{\Theta}(t),\label{zz1}\\
0&=&\dot{P}(t)x(t)+P(t)A(t)x(t)+P(t)B(t)u(t)\nonumber\\
&&+\hat{\Theta}(t)
+A'(t)P(t)x(t)+A'(t)\Theta(t)\nonumber\\
&&+\bar{A}'(t)q(t)+Q(t)x(t). \label{r6}
\end{eqnarray}
Using (\ref{zz1}) and (\ref{21}), (\ref{3}) becomes
\begin{eqnarray}
0&=&R(t)u(t)+B'(t)p(t)+\bar{B}'(t)q(t)\nonumber\\
&=&R(t)u(t)+B'(t)P(t)x(t)+B'(t)\Theta(t)\nonumber\\
&&+\bar{B}'(t)P(t)\bar{A}(t)x(t)+\bar{B}'(t)P(t)\bar{B}(t)u(t)\nonumber\\
&&+\bar{B}'(t)\bar{\Theta}(t)\nonumber\\
&=&\Upsilon_0(t)u(t)
+\Gamma_0(t)x(t)+B'(t)\Theta(t)\nonumber\\
&&+\bar{B}'(t)\bar{\Theta}(t),\label{25}
\end{eqnarray}
where $\Upsilon_0(t)$ and  $\Gamma_0(t)$ are respectively as in (\ref{Intr3}) and (\ref{Intr4}).
It should be noted that equilibrium condition (\ref{25}) is different from (\ref{Intr2}) when $\Theta(t)\not=0$.
We then have the following result.
\begin{theorem}\label{theorem1}
If it holds that
\begin{eqnarray}
Range \Big(\Gamma_0(t)\Big)\subseteq Range \Big(\Upsilon_0(t)\Big),\label{r23}
\end{eqnarray}
then Problem (Indefinite LQ) is solvable and the optimal solution
can be given by
\begin{eqnarray}
u(t)=-\Upsilon_0^{\dag}(t)\Gamma_0(t)x(t)+\Big(I-\Upsilon_0^{\dag}(t)\Upsilon_0(t)\Big)z(t),\label{r2}
\end{eqnarray}
where $z(t)\in R^m$ is arbitrary and $P(t)$ is given by the  Riccati equation (\ref{r1}):
In this case, $\Theta(t)=0$ holds in (\ref{21}).

\end{theorem}
\emph{Proof.}
If (\ref{r23}) holds, 
then $\Theta(t)=0, \bar{\Theta}(t)=0$. In fact,
substituting (\ref{zz1}) and (\ref{r2}) into (\ref{r6}) and using (\ref{r1}), we have
\begin{eqnarray}
\hat{\Theta}(t)&=&-\Big(A'(t)-\Gamma'_0(t)\Upsilon_0^{\dag}(t)B'(t)\Big)\Theta(t)\nonumber\\
&&-\Big(\bar{A}'(t)-\Gamma'_0(t)\Upsilon_0^{\dag}(t)\bar{B}'(t)\Big)\bar{\Theta}(t).\nonumber
\end{eqnarray}
Then it is obtained that
\begin{eqnarray}
d\Theta(t)&=&-\Big[\Big(A'(t)-\Gamma'_0(t)\Upsilon_0^{\dag}(t)B'(t)\Big)\Theta(t)+\Big(\bar{A}'(t)
\nonumber\\
&&-\Gamma'_0(t)\Upsilon_0^{\dag}(t)\bar{B}'(t)\Big)\bar{\Theta}(t)\Big]dt+\bar{\Theta}(t)dw(t).\nonumber
\end{eqnarray}
Using the terminal value $\Theta(T)=0,$ it follows that $\Theta(t)=0, \bar{\Theta}(t)=0.$
In this case, it is clear that $u(t)$ in (\ref{r2}) satisfies the equilibrium condition (\ref{25}) and thus (\ref{r2}) is one of the optimal controllers.
The proof is now completed.
\hfill $\blacksquare$

\begin{remark}
Conversely, it is easy to know that $ \Theta(t)\not=0$ if (\ref{r23}) does not hold, which is to be considered in the next section.
\end{remark}

%
%
%

\section{Solution to irregular LQ control}

This section focus on the stochastic LQ control problem when (\ref{r23}) does not hold, i.e.,
\begin{eqnarray}
Range \Big(\Gamma_0(t)\Big)\nsubseteq Range \Big(\Upsilon_0(t)\Big).\label{zz2}
\end{eqnarray}

In this case the problem is termed as the {\em irregular}  LQ control problem (IR-LQ). It is obvious that the controller can not be solved from the equilibrium condition (\ref{25}). However, as said in
Introduction,
this does not implies that the optimal control does not exists. This section aims to propose a multiple-layer optimization approach to solve the IR-LQ control of this case.


The proposed optimization approach in this paper may contain different layers.

\subsection{The First Layer}

From (\ref{zz2}), it is clear that $\Upsilon_0(t)$ is not invertible. We assume that $rank (\Upsilon_0(t))=m_0(t)< m$. Thus $rank\Big(I-\Upsilon_0^{\dag}(t)\Upsilon_0(t)\Big)=m-m_0(t)>0$. It is not difficult
to know that there is an elementary row transformation matrix $T_0(t)$ such that
\begin{eqnarray}
T_0(t)\Big(I-\Upsilon_0^{\dag}(t)\Upsilon_0(t)\Big)=\left[
                              \begin{array}{c}
                                 0  \\
                                 \Upsilon_{T_0}(t)\\
                              \end{array}
                            \right], \label{jnYC1}
\end{eqnarray}
where $\Upsilon_{T_0}(t)\in R^{[m-m_0(t)]\times m}$ is full row rank. Further denote
\begin{eqnarray}
 \left[ \begin{array}{cc} \ast & C_0'(t) \\
                              \end{array} \right] &=&\Gamma_0'(t)\Big(I-\Upsilon^{\dag}_0(t)\Upsilon_0(t)\Big){T_0}^{-1}(t),\nonumber\\
\left[ \begin{array}{cc}\ast & B_0(t)\\
                              \end{array}
                            \right] &=&B(t)\Big(I-\Upsilon_0^{\dag}(t)\Upsilon_0(t)\Big){T_0}^{-1}(t),\nonumber\\
\left[ \begin{array}{cc}\ast & \bar{B}_0(t)\\
                              \end{array}
                            \right] &=& \bar{B}(t)\Big(I-\Upsilon_0^{\dag}(t)\Upsilon_0(t)\Big) {T_0}^{-1}(t).\label{jn2}
\end{eqnarray}
Moreover, we also make the following denotations for the convenience of discussions
\begin{eqnarray}
A_0(t)&=&A(t)-B(t)\Upsilon_0^{\dag}(t)\Gamma_{0}(t),\nonumber\\
\bar{A}_0(t)&=&\bar{A}(t)-\bar{B}(t)\Upsilon_0^{\dag}(t)\Gamma_{0}(t),\nonumber\\
D_0(t)&=&-B(t)\Upsilon_0^{\dag}(t)B'(t),\nonumber\\
\bar{D}_0(t)&=&-\bar{B}(t)\Upsilon_0^{\dag}(t)B'(t),\nonumber\\
F_0(t)&=&-B(t)\Upsilon_0^{\dag}(t)\bar{B}'(t),\nonumber\\
\bar{F}_0(t)&=&-\bar{B}(t)\Upsilon_0^{\dag}(t)\bar{B}'(t).\nonumber
\end{eqnarray}

\begin{theorem}\label{theorem4}
Under the condition (\ref{zz2}), Problem (IR-LQ) is solvable if and only if there exists $u_1(t)\in R^{m-m_0(t)}$ such that
\begin{eqnarray}
0&=&C_0(t)x(t)+B_0'(t)\Theta(t)+\bar{B}_0'(t)\bar{\Theta}(t).\label{n2}
\end{eqnarray} holds,  where $u_1(t)$, $x(t),\Theta(t)$ and $\bar{\Theta}(t)$ satisfy the FBSDEs
\begin{eqnarray}
dx(t)&=&\Big(A_0(t)x(t)+D_0(t)\Theta(t)+F_0(t)\bar{\Theta}(t)\nonumber\\
&&+B_{0}(t)u_1(t)\Big)dt+\Big(\bar{A}_0(t)x(t)+\bar{D}_0(t)\Theta(t)\nonumber\\
&&+\bar{F}_0(t)\bar{\Theta}(t)+\bar{B}_0(t)u_1(t)\Big)dw(t),\label{n3}\\
d\Theta(t)&=&-\Big(A'_0(t)\Theta(t)+\bar{A}'_0(t)\bar{\Theta}(t)+C'_0(t)\nonumber\\
&&\times u_1(t)\Big)dt +\bar{\Theta}(t)dw(t), ~~~\Theta(T)=0 .\label{c8}
\end{eqnarray}

\end{theorem}

\emph{Proof.} ``Necessity".  If the problem is solvable, then from Lemma \ref{lemMP}, (\ref{25}) must hold. Thus, (\ref{25}) can be equivalently written as
\begin{eqnarray}
u(t)&=&-\Upsilon_0^{\dag}(t)\Big(\Gamma_0(t)x(t)+B'(t)\Theta(t)+\bar{B}'(t)\bar{\Theta}(t)\Big)\nonumber\\
&&+\Big(I-\Upsilon_0^{\dag}(t)\Upsilon_0(t)\Big)z(t),\label{n1}
\end{eqnarray}
where $z(t)$ is a vector with compatible dimension such that following equality hold
\begin{eqnarray}
0&=&\Big(I-\Upsilon_0(t)\Upsilon_0^{\dag}(t)\Big)\Big(\Gamma_0(t)x(t)+B'(t)\Theta(t)\nonumber\\
&&+\bar{B}'(t)\bar{\Theta}(t)\Big). \label{nz1}
\end{eqnarray}
Denote \begin{eqnarray}
T_0(t)\Big(I-\Upsilon_0^{\dag}(t)\Upsilon_0(t)\Big)z(t)=\left[
                              \begin{array}{c}
                                 0  \\
                                 u_1(t)\\
                              \end{array}
                            \right], \label{jn1}
\end{eqnarray}
where $u_1(t)=\Upsilon_{T_0}(t) z(t)\in R^{m-m_0(t)}$. Note that $\Upsilon_{T_0}(t)\in R^{[m-m_0(t)]\times m}$ is full row rank.
Now we rewrite (\ref{nz1}) as (\ref{n2}). First, note that
\begin{eqnarray}
&&I-\Upsilon_0(t)\Upsilon_0^{\dag}(t)\nonumber\\
&=&\Big(I-\Upsilon_0(t)\Upsilon_0^{\dag}(t)\Big)\Big(I-\Upsilon_0(t)\Upsilon_0^{\dag}(t)\Big)\nonumber \\
&=&\Big(I-\Upsilon_0(t)\Upsilon_0^{\dag}(t)\Big)T'_0(t)\Big(T^{-1}_0(t)\Big)'\Big(I-\Upsilon_0(t)\Upsilon_0^{\dag}(t)\Big)\nonumber \\
&=&\left[
\begin{array}{cc}
                                 0 &
                                 \Upsilon'_{T_0}(t)\\
                              \end{array}
                            \right]\Big(T^{-1}_0(t)\Big)'\Big(I-\Upsilon_0(t)\Upsilon_0^{\dag}(t)\Big),
\end{eqnarray}
where (\ref{jnYC1}) has been used in the derivation of the last equality. By using (\ref{jn2}), (\ref{nz1}) can be written as
\begin{eqnarray}
0&=&\Upsilon'_{T_0}(t)\Big[C_0(t)x(t)+B_0'(t)\Theta(t)+\bar{B}_0'(t)\bar{\Theta}(t)\Big].\label{nYc2}
\end{eqnarray}
Note that $\Upsilon'_{T_0}(t)$ is full column rank, (\ref{nYc2}) is rewritten as (\ref{n2}) directly.
By substituting (\ref{n1}) and (\ref{zz1}) into (\ref{r6}) and using (\ref{r1}), it yields that
\begin{eqnarray}
0&=&\dot{P}(t)x(t)+P(t)A(t)x(t)+\hat{\Theta}(t)
+A'(t)P(t)x(t)\nonumber\\
&&+A'(t)\Theta(t)+\bar{A}'(t)P(t)\bar{A}(t)x(t)+\bar{A}'(t)\bar{\Theta}(t)\nonumber\\
&&+Q(t)x(t)-\Gamma_0'(t)\Upsilon_0^{\dag}(t)\Big(\Gamma_0(t)x(t)+B'(t)\Theta(t)\nonumber\\
&&+\bar{B}'(t)\bar{\Theta}(t)\Big)+\Gamma_0'(t)\Big(I-\Upsilon_0^{\dag}(t)\Upsilon_0(t)\Big)z(t)\nonumber\\
&=&\hat{\Theta}(t)+\Big(A'(t)-\Gamma_0'(t)\Upsilon_0^{\dag}(t)B'(t)\Big)\Theta(t)\nonumber\\
&&+\Big(\bar{A}'(t)-\Gamma_0'(t)\Upsilon_0^{\dag}(t)\bar{B}'(t)\Big)\bar{\Theta}(t)\nonumber\\
&&+\Gamma_0'(t)\Big(I-\Upsilon_0^{\dag}(t)\Upsilon_0(t)\Big)z(t).\label{ZJ1}
\end{eqnarray}
In view of the fact that $\Big(I-\Upsilon_0^{\dag}(t)\Upsilon_0(t)\Big)^2=I-\Upsilon_0^{\dag}(t)\Upsilon_0(t),$ it is obtained that
\begin{eqnarray}
&&\Gamma_0'(t)\Big(I-\Upsilon_0^{\dag}(t)\Upsilon_0(t)\Big)z(t)\nonumber\\
&=&\Gamma_0'(t)\Big(I-\Upsilon_0^{\dag}(t)\Upsilon_0(t)\Big)T_0^{-1}(t)T_0(t)\Big(I-\Upsilon_0^{\dag}(t)\nonumber\\
&&\times\Upsilon_0(t)\Big)z(t)\nonumber\\
&=&\Gamma_0'(t)\Big(I-\Upsilon_0^{\dag}(t)\Upsilon_0(t)\Big)T_0^{-1}(t)\left[
                              \begin{array}{c}
                                 0  \\
                                 u_1(t)\\
                              \end{array}
                            \right]\nonumber\\
&=&\left[
     \begin{array}{cc}
       * & C_0'(t) \\
     \end{array}
   \right]\left[
                              \begin{array}{c}
                                 0  \\
                                 u_1(t)\\
                              \end{array}
                            \right]\nonumber\\
&=&C_0'(t)u_1(t).\label{YCJ1}
\end{eqnarray}

Thus, from ({\ref{ZJ1}) we have
\begin{eqnarray}
\hat{\Theta}(t)=-\Big[A'_0(t)\Theta(t)+\bar{A}'_0(t)\bar{\Theta}(t)+C'_0(t) u_1(t)\Big],
\end{eqnarray}
this implies that the dynamic of $\Theta$ is as (\ref{c8}) using (\ref{zj1}). By substituting (\ref{n1}) into (\ref{1}),
one has the dynamic (\ref{n3}) of the state.

``Sufficiency" We now show Problem (IR-LQ) is solvable if there exists $u_1(t)$ to achieve (\ref{n2}). In fact, if (\ref{n2}) is true then (\ref{n1}) and (\ref{nz1}) can be rewritten as (\ref{25}).  Further, by taking reverse procedures to (\ref{21})-(\ref{25}), it is easily verified that
$p(t)=P(t)x(t)+\Theta(t),q(t)=P(t)\big[\bar{A}(t)x(t)+\bar{B}(t)u(t)\big]+\bar{\Theta}(t)$, where $x, \Theta,\bar{\Theta}$
satisfy (\ref{n2})-(\ref{c8}), are the solutions to (\ref{3})-(\ref{2}). Thus, Problem (IR-LQ) is solvable according to
Lemma \ref{lemMP}.
The proof is now completed. \hfill $\blacksquare$

\begin{remark}
Applying (\ref{jn1}) we have
\begin{eqnarray}
\Big(I-\Upsilon_0^{\dag}(t)\Upsilon_0(t)\Big)z(t)&=&T^{-1}_0(t)\left[
                              \begin{array}{c}
                                 0  \\
                                 u_1(t)\\
                              \end{array}
                            \right]\nonumber\\
 &=&G_0(t) u_1(t), \label{jny1}
\end{eqnarray}
where $
 \left[ \begin{array}{cc} \ast & G_0(t) \\
                              \end{array} \right] ={T_0}^{-1}(t)$.  Thus (\ref{n1}) can be rewritten as

\begin{eqnarray}
u(t)&=&-\Upsilon_0^{\dag}(t)\Big(\Gamma_0(t)x(t)+B'(t)\Theta(t)+\bar{B}'(t)\bar{\Theta}(t)\Big)\nonumber\\
&&+G_0(t)u_1(t). \label{rY2}
\end{eqnarray}

\end{remark}

\begin{observation}\label{obser1}
The solution to FBSDEs (\ref{3})-(\ref{2}) is homogeneous, i.e., $\Theta(t) =0,$ if and only if the Riccati equation (\ref{r1}) is regular.
\end{observation}
\emph{Proof.} ``Necessity" Assume that the Riccati equation (\ref{r1}) is not regular, then
\begin{eqnarray}
\Big(I-\Upsilon_0(t)\Upsilon_0^{\dag}(t)\Big)\Gamma_0(t)\not=0, \label{qd3}
\end{eqnarray}
which implies $C_0(t)\not=0$ using (\ref{YCJ1}). Thus according to
Theorem \ref{theorem4}, $\Theta(t)$ satisfies (\ref{c8}) which is not equal to zero. This is a contradiction.

``Sufficiency" The necessity follows from Theorem \ref{theorem1}.
The proof is now completed. \hfill $\blacksquare$
\\

Now Problem IR-LQ is converted into finding $u_1(t)$ to achieve (\ref{n2}) associated with (\ref{n3})-(\ref{c8}). To solve the problem, without loss of generality, we let
\begin{eqnarray}
\Theta_1(t)=\Theta(t)-P_1(t)x(t),\label{c48}
\end{eqnarray}
where $P_1(t)$ obeys
\begin{eqnarray}
0&=&\dot{P}_1(t)+P_1(t)A_0(t)+A_0'(t)P_1(t)+P_1(t)D_0(t)\nonumber\\
&&\times P_1(t)+\Big(\bar{A}'_0(t)+P_1(t)F_0(t)\Big)\Big(I-P_1(t)\bar{F}_{0}(t)\Big)^{\dag}\nonumber\\
&&\times P_1(t)\Big(\bar{A}_0(t)
+\bar{D}_0(t)P_1(t)\Big)-\Gamma_1'(t)\Upsilon_1^{\dag}(t)\Gamma_1(t), \nonumber\\\label{p1}
\end{eqnarray}
where
\begin{eqnarray}
\Upsilon_1(t)&=&\bar{B}_{0}'(t)\Big(I-P_1(t)\bar{F}_{0}(t)\Big)^{\dag}P_1(t)\bar{B}_{0}(t), \label{pz1}\\
\Gamma_1(t)&=&C_0(t)+B_{0}'(t)P_1(t)+\bar{B}_{0}'(t)\Big(I-P_1(t)\nonumber\\
&&\times \bar{F}_{0}(t)\Big)^{\dag}P_1(t)\Big(\bar{A}_0(t)+\bar{D}_{0}(t) P_1(t)\Big),\label{pz2}
\end{eqnarray}
with arbitrary terminal value $P_1(T)$.

Similarly, we assume that
\begin{eqnarray}
d\Theta_1(t)=\hat{\Theta}_1(t)dt+\bar{\Theta}_1(t)dw(t),
\end{eqnarray}
where $\hat{\Theta}_1(t)$ and $\bar{\Theta}_1(t)$ are to be determined. Then by taking It\^{o}'s formula to (\ref{c48}), it is obtained that
\begin{eqnarray}
d\Theta(t)
&=&\dot{P}_1(t)x(t)dt+P_1(t)dx(t)+d\Theta_1(t)\nonumber\\
&=&\dot{P}_1(t)x(t)dt+P_1(t)\Big[A_0(t)x(t)
+D_0(t)\Theta(t)\nonumber\\
&&+F_0(t)\bar{\Theta}(t)+B_{0}(t)u_1(t)\Big]dt+P_1(t)\Big[\bar{A}_0(t)\nonumber\\
&&\times x(t)
+\bar{D}_0(t)\Theta(t)+\bar{F}_0(t)\bar{\Theta}(t)+\bar{B}_{0}(t)\nonumber\\
&&\times u_1(t)\Big]dw(t)+\hat{\Theta}_1(t)dt+\bar{\Theta}_1(t)dw(t).
\end{eqnarray}
Combining with (\ref{c8}), we have
\begin{eqnarray}
0&=&\dot{P}_1(t)x(t)+P_1(t)\Big[A_0(t)x(t)
+D_0(t)\Theta(t)\nonumber\\
&&+F_0(t)\bar{\Theta}(t)+B_{0}(t)u_1(t)\Big]+A'_0(t)P_1(t)x(t)\nonumber\\
&&+\hat{\Theta}_1(t)+A'_0(t)\Theta_1(t)+\bar{A}'_0(t)\bar{\Theta}(t)\nonumber\\
&&+C_0'(t)u_1(t),\label{c49}\\
\bar{\Theta}(t)&=&\bar{\Theta}_1(t)+P_1(t)\Big[\bar{A}_0(t)x(t)
+\bar{D}_0(t)\Theta(t)\nonumber\\
&&+\bar{F}_0(t)\bar{\Theta}(t)+\bar{B}_{0}(t)u_1(t)\Big].\label{c50}
\end{eqnarray}
Thus it is obtained from (\ref{c50}) that
\begin{eqnarray}
\Big(I-P_1(t)\bar{F}_{0}(t)\Big)\bar{\Theta}(t)
&=&\Big[\bar{\Theta}_1(t)+P_1(t)\Big(\bar{A}_0(t)x(t)\nonumber\\
&&
+\bar{D}_0(t)\Theta(t)+\bar{B}_{0}(t)u_1(t)\Big)\Big]. \nonumber
\end{eqnarray}
Note that $\bar{\Theta}(t)$, satisfying (\ref{c8}), is unique, which is given from the above equation as
\begin{eqnarray}
\bar{\Theta}(t)&=&\Big(I-P_1(t)\bar{F}_{0}(t)\Big)^{\dag}\Big[\bar{\Theta}_1(t)+P_1(t)\Big(\bar{A}_0(t)x(t)\nonumber\\
&&+\bar{D}_0(t)\Theta(t)+\bar{B}_{0}(t)u_1(t)\Big)\Big]+\Big[I-\Big(I-P_1(t)\nonumber\\
&&\times\bar{F}_{0}(t)\Big)^{\dag}\Big(I-P_1(t)\bar{F}_{0}(t)\Big)\Big]\varphi(t),\label{m26}
\end{eqnarray}
where $\varphi(t)$ is a vector with compatible dimension. Substituting (\ref{m26}) into (\ref{25}), it yields that
\begin{eqnarray}
0
&=&\Upsilon_0(t)u(t)
+\Gamma_0(t)x(t)+B'(t)P_1(t)x(t)+B'(t)\nonumber\\
&&\times\Theta_1(t)+\bar{B}'(t)\Big(I-P_1(t)\bar{F}_{0}(t)\Big)^{\dag}\Big[\bar{\Theta}_1(t)\nonumber\\
&&+P_1(t)\Big(\bar{A}_0(t)x(t)+\bar{D}_0(t)\Theta(t)+\bar{B}_{0}(t)u_1(t)\Big)\Big]\nonumber\\
&&+\bar{B}'(t)\Big[I-\Big(I-P_1(t)\bar{F}_{0}(t)\Big)^{\dag}\Big(I-P_1(t)\nonumber\\
&&\times\bar{F}_{0}(t)\Big)\Big]\varphi(t). \label{c61}
\end{eqnarray}
It is not difficult to know that following condition holds
\begin{eqnarray}
0&=&\bar{B}'(t)\Big[I-\Big(I-P_1(t)\bar{F}_{0}(t)\Big)^{\dag}\Big(I-P_1(t)\bar{F}_{0}(t)\Big)\Big].\nonumber\\\label{z17}
\end{eqnarray}
In fact, if (\ref{z17}) does not hold, then there exists an unknown vector $\varphi(t)$ in the equilibrium condition (\ref{c61}),
that is, an unknown vector is involved in the solution to $u(t).$ This implies that Problem (IR-LQ) is unsolvable.

\begin{remark}\label{rem3}
Using (\ref{z17}) and the definitions of $\bar{B}_0(t), \bar{D}_0(t), F_0(t), \bar{F}_0(t)$, it is immediate to obtain that
\begin{eqnarray}
0=L(t)\Big[I-\Big(I-P_1(t)\bar{F}_{0}(t)\Big)^{\dag}\Big(I-P_1(t)\bar{F}_{0}(t)\Big)\Big],\label{z1}
\end{eqnarray}
where $L(t)$ may be $\bar{B}'_0(t), \bar{D}'_0(t), F_0(t), \bar{F}_0(t)$.

\end{remark}

Substituting (\ref{m26}) into (\ref{n2})  yields that
\begin{eqnarray}
0&=&C_0(t)x(t)+B_{0}'(t)\Theta(t)+\bar{B}_{0}'(t)\bar{\Theta}(t)\nonumber\\
&=& \Big[\bar{B}_{0}'(t)\Big(I-P_1(t)\bar{F}_{0}(t)\Big)^{\dag} P_1(t)\bar{B}_{0}(t)\Big]u_1(t)\nonumber\\
&&+\Big[C_{0}(t)+B_{0}'(t)P_1(t)\nonumber\\
&&+\bar{B}_{0}'(t)\Big(I-P_1(t)\bar{F}_{0}(t)\Big)^{\dag}P_1(t)\Big(\bar{A}_0(t)+\bar{D}_{0}(t)\nonumber\\
&&\times P_1(t)\Big)\Big]x(t)+H_{1}'(t)\Theta_1(t)+\bar{H}_{1}'(t)\bar{\Theta}_1(t)\nonumber\\
&&
+ \bar{B}_0'(t)\Big[I-\Big(I-P_1(t)\bar{F}_{0}(t)\Big)^{\dag}\Big(I-P_1(t)\nonumber\\
&&\times\bar{F}_{0}(t)\Big)\Big]\varphi(t)\nonumber\\
&=&\Upsilon_1(t)u_1(t)+\Gamma_1(t)x(t)+H_{1}'(t)\Theta_1(t)\nonumber\\
&&+\bar{H}_{1}'(t)\bar{\Theta}_1(t)+\bar{B}_0'(t)\Big[I-\Big(I-P_1(t)\bar{F}_{0}(t)\Big)^{\dag}\nonumber\\
&&\times\Big(I-P_1(t)\bar{F}_{0}(t)\Big)\Big]\varphi(t)\nonumber\\
&=&\Upsilon_1(t)u_1(t)+\Gamma_1(t)x(t)+H_1'(t)\Theta_1(t)\nonumber\\
&&+\bar{H}_{1}'(t)\bar{\Theta}_1(t),\label{c52}
\end{eqnarray}
where the matrices $\Upsilon_1(t)$ and $\Gamma_1(t)$ are respectively as (\ref{pz1}) and (\ref{pz2}), and
\begin{eqnarray}
H_{1}'(t)&=&B_{0}'(t)+\bar{B}_{0}'(t)\Big(I-P_1(t)\bar{F}_0(t)\Big)^{\dag}P_1(t)\bar{D}_0(t),\nonumber\\
\bar{H}_{1}'(t)&=&\bar{B}_{0}'(t)\Big(I-P_1(t)\bar{F}_0(t)\Big)^{\dag}.\nonumber
\end{eqnarray}

In the sequel, there are two cases to be considered, one is the case of $Range \Big(\Gamma_1(t)\Big)\subseteq Range \Big(\Upsilon_1(t)\Big)$
and the other is $Range \Big(\Gamma_1(t)\Big)\nsubseteq Range \Big(\Upsilon_1(t)\Big).$

We firstly consider the case that $Range \Big(\Gamma_1(t)\Big)\subseteq Range \Big(\Upsilon_1(t)\Big).$ In this case,
\begin{eqnarray}
\Gamma_1'(t)\Big(I-\Upsilon_1^{\dag}(t)\Upsilon_1(t)\Big)=0,\label{ga1}
\end{eqnarray}
and $u_1(t)$ is obtained from  (\ref{c52}) as,
\begin{eqnarray}
u_1(t)&=&-\Upsilon_1^{\dag}(t)\Big(\Gamma_1(t)x(t)+H_{1}'(t)\Theta_1(t)+\bar{H}_{1}'(t)\bar{\Theta}_1(t)\Big)\nonumber\\
&&+\Big(I-\Upsilon_1^{\dag}(t)\Upsilon_1(t)\Big)z_1(t),\label{m24}
\end{eqnarray}
where $z_1(t)$ is an arbitrary vector with compatible dimension.
By substituting (\ref{c48}), (\ref{m26}) and (\ref{m24}) into the dynamic of the system (\ref{n3}), we have
\begin{eqnarray}
dx(t)&=&\Big\{\Big(M_1(t)-\hat{H}_{1}(t)\Upsilon_1^{\dag}(t)\Gamma_1(t)\Big)x(t)\nonumber\\
&&+\Big[D_0(t)+F_0(t)\Big(I-P_1(t)\bar{F}_0(t)\Big)^{\dag}P_1(t)\bar{D}_{0}(t)\nonumber\\
&&-\hat{H}_{1}(t)\Upsilon_1^{\dag}(t)H_{1}'(t)\Big]\Theta_1(t)
+\Big[F_0(t)\Big(I-P_1(t)\nonumber\\
&&\times\bar{F}_0(t)\Big)^{\dag}-\hat{H}_{1}(t)\Upsilon_1^{\dag}(t)\bar{H}_{1}'(t)\Big]\bar{\Theta}_1(t)\nonumber\\
&&+\hat{H}_{1}(t)\Big(I-\Upsilon_1^{\dag}(t)\Upsilon_1(t)\Big)z_1(t)\Big\}dt+\Big\{\Big(\bar{M}_1(t)\nonumber\\
&&-\hat{\bar{H}}_{1}(t)\Upsilon_1^{\dag}(t)\Gamma_1(t)\Big)x(t)+\hat{\bar{H}}_{1}(t) \Big(I-\Upsilon_1^{\dag}(t)\nonumber\\
&&\times\Upsilon_1(t)\Big)z_1(t)+\Big[\bar{D}_0(t)+\bar{F}_0(t)\Big(I-P_1(t)\nonumber\\
&&\times \bar{F}_0(t)\Big)^{\dag}P_1(t)\bar{D}_{0}(t)-\hat{\bar{H}}_{1}(t)\Upsilon_1^{\dag}(t)H_{1}'(t)\Big]\Theta_1(t)\nonumber\\
&&
+\Big[\bar{F}_0(t)\Big(I-P_1(t)\bar{F}_0(t)\Big)^{\dag}\nonumber\\
&&-\hat{\bar{H}}_{1}(t)\Upsilon_1^{\dag}(t)\bar{H}_{1}'(t)\Big]\bar{\Theta}_1(t)\Big\}dw(t),\label{m25}
\end{eqnarray}
where
\begin{eqnarray}
M_1(t)&=&A_0(t)
+D_{0}(t)P_1(t)+F_0(t)[I-P_1(t)\bar{F}_{0}(t)]^{\dag} \nonumber\\
&&\times P_1(t)[\bar{A}_0(t)+\bar{D}_0(t)P_1(t)]\nonumber\\
\bar{M}_1(t)&=&\bar{A}_0(t)
+\bar{D}_{0}(t)P_1(t)+\bar{F}_0(t)[I-P_1(t)\bar{F}_{0}(t)]^{\dag}\nonumber\\
&&\times P_1(t)[\bar{A}_0(t)+\bar{D}_0(t)P_1(t)]\nonumber\\
\hat{H}_{1}(t)&=&B_{0}(t)+F_0(t)\Big(I-P_1(t)\bar{F}_0(t)\Big)^{\dag}P_1(t)\bar{B}_{0}(t),\nonumber\\
\hat{\bar{H}}_{1}(t)&=&\bar{B}_{0}(t)+\bar{F}_0(t)\Big(I-P_1(t)\bar{F}_0(t)\Big)^{\dag}P_1(t)\bar{B}_{0}(t).\nonumber
\end{eqnarray}
Similarly, substituting (\ref{c48}), (\ref{m26}) into  (\ref{c49}), we have that
\begin{eqnarray}
0&=&\Big[\dot{P}_1(t)+P_1(t)A_0(t)+A'_0(t)P_1(t)
+P_1(t)D_0(t)\nonumber\\
&&\times P_1(t)+\Big(\bar{A}'_0(t)+P_1(t)F_0(t)\Big)\Big(I-P_1(t)\bar{F}_{0}(t)\Big)^{\dag}\nonumber\\
&&\times P_1(t)\Big(\bar{A}_0(t)+\bar{D}_0(t)P_1(t)\Big)\Big]x(t)+\Big[A'_0(t)\nonumber\\
&&+P_1(t)D_0(t)+\Big(\bar{A}'_0(t)+P_1(t)F_0(t)\Big)\Big(I-P_1(t)\nonumber\\
&&\times \bar{F}_{0}(t)\Big)^{\dag}P_1(t)\bar{D}_0(t)\Big]\Theta_1(t)+\Big[C_0'(t)+P_1(t)B_{0}(t)\nonumber\\
&&+\Big(\bar{A}'_0(t)+P_1(t)F_0(t)\Big)\Big(I-P_1(t)\bar{F}_{0}(t)\Big)^{\dag}P_1(t)\nonumber\\
&&\times \bar{B}_{0}(t)\Big]u_1(t)+\hat{\Theta}_1(t)+\Big(\bar{A}'_0(t)+P_1(t)F_0(t)\Big)\nonumber\\
&&\times \Big(I-P_1(t)\bar{F}_{0}(t)\Big)^{\dag}\bar{\Theta}_1(t) +\Big(\bar{A}'_0(t)+P_1(t)F_0(t)\Big)\nonumber\\
&&\times \Big[I-\Big(I-P_1(t)\bar{F}_{0}(t)\Big)^{\dag}\Big(I-P_1(t)\bar{F}_{0}(t)\Big)\Big]\varphi(t).\nonumber\\\label{c60}
\end{eqnarray}
By denoting
\begin{eqnarray}
\hat{M}_1(t)&=&A_0(t)
+D_{0}(t)P_1(t)+F_0(t)P_1(t)[I-\bar{F}_{0}(t)\nonumber\\
&&\times P_1(t)]^{\dag} [\bar{A}_0(t)+\bar{D}_0(t)P_1(t)]\nonumber\\
\hat{\bar{M}}_1(t)&=&[I-\bar{F}_{0}(t)P_1(t)]^{\dag} [\bar{A}_0(t)+\bar{D}_0(t)P_1(t)]\nonumber\\
\hat{\Gamma}_{1}(t)&=&C_0(t)+B_{0}'(t)P_1(t)+\bar{B}_{0}'(t)P_1(t)\Big(I-\nonumber\\
&& \bar{F}_{0}(t)P_1(t)\Big)^{\dag}\Big(\bar{A}_0(t)+\bar{D}_{0}(t) P_1(t)\Big)
\end{eqnarray}
(\ref{c60}) is reduced to
\begin{eqnarray}
0&=&\Big[\dot{P}_1(t)+P_1(t)A_0(t)+A'_0(t)P_1(t)
+P_1(t)D_0(t)\nonumber\\
&&\times P_1(t)+\Big(\bar{A}'_0(t)+P_1(t)F_0(t)\Big)\Big(I-P_1(t)\bar{F}_{0}(t)\Big)^{\dag}\nonumber\\
&&\times P_1(t)\Big(\bar{A}_0(t)+\bar{D}_0(t)P_1(t)\Big)\Big]x(t)+\hat{M}_1'(t)\Theta_1(t)\nonumber\\
&&+\hat{\Gamma}_1'(t)u_1(t)+\hat{\Theta}_1(t)+\hat{\bar{M}}_1'(t)\bar{\Theta}_1(t)+\bar{A}'_0(t)\nonumber\\
&&\times\Big[I-\Big(I-P_1(t)\bar{F}_{0}(t)\Big)^{\dag}\Big(I-P_1(t)\bar{F}_{0}(t)\Big)\Big]\varphi(t),\nonumber
\end{eqnarray}
By substituting (\ref{m24}) into the above equation, we further obtain that
\begin{eqnarray}
0&=&\Big[\dot{P}_1(t)+P_1(t)A_0(t)+A'_0(t)P_1(t)
+P_1(t)D_0(t)\nonumber\\
&&\times P_1(t)+\Big(\bar{A}'_0(t)+P_1(t)F_0(t)\Big)\Big(I-P_1(t)\bar{F}_{0}(t)\Big)^{\dag}\nonumber\\
&&\times P_1(t)\Big(\bar{A}_0(t)+\bar{D}_0(t)P_1(t)\Big)-\hat{\Gamma}_1'(t)\Upsilon_1^{\dag}(t)\Gamma_1(t)\Big]\nonumber\\
&&\times x(t)+\Big(\hat{M}_1'(t)-\Gamma_1'(t)\Upsilon_1^{\dag}(t)H_{1}(t)\Big)\Theta_1(t)\nonumber\\
&&
+\Big(\hat{\bar{M}}_1'(t)-\Gamma_1'(t)\Upsilon_1^{\dag}(t)\bar{H}_{1}(t)\Big)\bar{\Theta}_1(t)\nonumber\\
&&+\hat{\Gamma}_1'(t)\Big(I-\Upsilon_1^{\dag}(t)\Upsilon_1(t)\Big)z_1(t)+\hat{\Theta}_1(t)+\bar{A}'_0(t)\nonumber\\
&&\times\Big[I-\Big(I-P_1(t)\bar{F}_{0}(t)\Big)^{\dag}\Big(I-P_1(t)\bar{F}_{0}(t)\Big)\Big]\varphi(t).\nonumber\\\label{au2}
\end{eqnarray}
 Similar to the derivation of (\ref{z17}), the following equality is necessary for the solvability of optimal control.
\begin{eqnarray}
0&=&\bar{A}'_0(t)\Big[I-\Big(I-P_1(t)\bar{F}_{0}(t)\Big)^{\dag}\Big(I-P_1(t)\bar{F}_{0}(t)\Big)\Big].\nonumber\\\label{z3}
\end{eqnarray}

\begin{lemma}\label{leminv}
Under (\ref{z17}) and (\ref{z3}), it holds that
\begin{enumerate}
  \item  Commutative law
  \begin{eqnarray}
  &&L_1'(t)P_1(t)\Big(I-\bar{F}_0(t)P_1(t)\Big)^{\dag}L_2(t)\nonumber\\
  &=&L_{1}'(t)\Big(I-P_1(t)\bar{F}_0(t)\Big)^{\dag}P_1(t)L_{2}(t),\nonumber
  \end{eqnarray}
  where $L_1(t),L_2(t)$ may be $\bar{B}_{0}(t),\bar{D}_0(t),\bar{A}_0(t),\bar{F}_0(t).$
  \item Formula of More-Penrose inverse for sum of matrices
  \begin{eqnarray}
  &&\Big(I-\bar{F}_0(t)P_1(t)\Big)^{\dag}L(t)\nonumber\\
  &=&\Big[I+\bar{F}_0(t)\Big(I-P_1(t)\bar{F}_{0}(t)\Big)^{\dag}P_1(t)\Big]L(t),\nonumber
  \end{eqnarray} where $L(t)$ may be $\bar{B}_{0}(t),\bar{D}_0(t),\bar{A}_0(t),\bar{F}_0(t).$

\end{enumerate}
\end{lemma}
\emph{Proof.}
\begin{enumerate}
  \item We firstly consider the case that $L_1(t)=\bar{D}_{0}(t)$ and $L_2(t)=\bar{B}_{0}(t).$
  By using Remark \ref{rem3}, we have
  \begin{eqnarray}
  &&\bar{D}_{0}'(t)P_1(t)\Big(I-\bar{F}_0(t)P_1(t)\Big)^{\dag}\bar{B}_{0}(t)\nonumber\\
  &=&\bar{D}_{0}'(t)\Big(I-P_1(t)\bar{F}_{0}(t)\Big)^{\dag}\Big(I-P_1(t)\bar{F}_{0}(t)\Big)P_1(t)\nonumber\\
  &&\times \Big(I-\bar{F}_0(t)P_1(t)\Big)^{\dag}\bar{B}_{0}(t)\nonumber\\
  &=&\bar{D}_{0}'(t)\Big(I-P_1(t)\bar{F}_{0}(t)\Big)^{\dag}P_1(t)\Big(I-\bar{F}_{0}(t)P_1(t)\Big)\nonumber\\
  &&\times \Big(I-\bar{F}_0(t)P_1(t)\Big)^{\dag}\bar{B}_{0}(t)\nonumber\\
  &=&\bar{D}_{0}'(t)\Big(I-P_1(t)\bar{F}_{0}(t)\Big)^{\dag}P_1(t)\bar{B}_{0}(t),\nonumber
  \end{eqnarray}
  where (\ref{z1}) has been used in the derivation of the last equality. The other cases can be obtained similarly.
  \item Consider the case that $L(t)=\bar{B}_{0}(t).$ By Remark \ref{rem3},  we have
  \begin{eqnarray}
&&\Big(I-\bar{F}_0(t)P_1(t)\Big)^{\dag}\bar{B}_{0}(t)\nonumber\\
&=&\bar{B}_{0}(t)+\bar{F}_0(t)P_1(t)\Big(I-\bar{F}_0(t)P_1(t)\Big)^{\dag}\bar{B}_{0}(t)\nonumber\\
&=&\bar{B}_{0}(t)+\bar{F}_0(t)\Big(I-P_1(t)\bar{F}_{0}(t)\Big)^{\dag}\Big(I-P_1(t)\nonumber\\
&&\times \bar{F}_{0}(t)\Big)P_1(t)\Big(I-\bar{F}_0(t)P_1(t)\Big)^{\dag}\bar{B}_{0}(t)\nonumber\\
&=&\bar{B}_{0}(t)+\bar{F}_0(t)\Big(I-P_1(t)\bar{F}_{0}(t)\Big)^{\dag}P_1(t)\Big(I-\nonumber\\
&&\bar{F}_{0}(t)P_1(t)\Big)\Big(I-\bar{F}_0(t)P_1(t)\Big)^{\dag}\bar{B}_{0}(t)\nonumber\\
&=&\bar{B}_{0}(t)+\bar{F}_0(t)\Big(I-P_1(t)\bar{F}_{0}(t)\Big)^{\dag}P_1(t)\bar{B}_{0}(t),\nonumber
  \end{eqnarray}
where (\ref{z1}) has been used in the derivation of the last equality. The derivations of the cases that $L(t)=\bar{D}_0(t),\bar{A}_0(t),\bar{F}_0(t)$ are
similar. So we omit. The proof is now completed. \hfill $\blacksquare$

\end{enumerate}

\begin{remark}\label{rem1}
Based on Lemma \ref{leminv}, we have that
\begin{eqnarray}
\hat{H}_{1}(t)
&=&B_{1}(t)\nonumber\\
\hat{\bar{H}}_{1}(t)
&=&\bar{H}_{1}(t),\nonumber\\
\hat{M}_1(t)&=&M_1(t),\nonumber\\
\hat{\bar{M}}_1(t)&=&\bar{M}_1(t),\nonumber\\
\hat{\Gamma}_1(t)&=&\Gamma_1(t),\nonumber\\
\bar{M}_1'(t)&=&[\bar{A}_0(t)+\bar{D}_0(t)P_1(t)]'[I-P_1(t)\bar{F}_{0}(t)]^{\dag},\nonumber\\
\bar{H}_{1}'(t)
&=&\bar{B}_{0}'(t)+\bar{B}_{0}'(t)\Big(I-P_1(t)\bar{F}_{0}(t)\Big)^{\dag}P_1(t)\bar{F}_0(t).\nonumber
\end{eqnarray}

\end{remark}

Based on Remark \ref{rem1}, we rewrite (\ref{au2}) as
\begin{eqnarray}
0&=&\Big[\dot{P}_1(t)+P_1(t)A_0(t)+A'_0(t)P_1(t)
+P_1(t)D_0(t)\nonumber\\
&&\times P_1(t)+\Big(\bar{A}'_0(t)+P_1(t)F_0(t)\Big)\Big(I-P_1(t)\bar{F}_{0}(t)\Big)^{\dag}\nonumber\\
&&\times P_1(t)\Big(\bar{A}_0(t)+\bar{D}_0(t)P_1(t)\Big)-\Gamma_1'(t)\Upsilon_1^{\dag}(t)\Gamma_1(t)\Big]\nonumber\\
&&\times x(t)+\Big(M_1'(t)-\Gamma_1'(t)\Upsilon_1^{\dag}(t)H_{1}(t)\Big)\Theta_1(t)\nonumber\\
&&
+\Big(\bar{M}_1'(t)-\Gamma_1'(t)\Upsilon_1^{\dag}(t)\bar{H}_{1}(t)\Big)\bar{\Theta}_1(t)\nonumber\\
&&+\Gamma_1'(t)\Big(I-\Upsilon_1^{\dag}(t)\Upsilon_1(t)\Big)z_1(t)+\hat{\Theta}_1(t)+\bar{A}'_0(t)\nonumber\\
&&\times\Big[I-\Big(I-P_1(t)\bar{F}_{0}(t)\Big)^{\dag}\Big(I-P_1(t)\bar{F}_{0}(t)\Big)\Big]\varphi(t)\nonumber \\
&=&\Big(M_1'(t)-\Gamma_1'(t)\Upsilon_1^{\dag}(t)H_{1}(t)\Big)\Theta_1(t)\nonumber\\
&&
+\Big(\bar{M}_1'(t)-\Gamma_1'(t)\Upsilon_1^{\dag}(t)\bar{H}_{1}(t)\Big)\bar{\Theta}_1(t)\nonumber\\
&&+\Gamma_1'(t)\Big(I-\Upsilon_1^{\dag}(t)\Upsilon_1(t)\Big)z_1(t)+\hat{\Theta}_1(t),\label{au3}
\end{eqnarray}
where (\ref{p1}) has been inserted in the last equality. Together with (\ref{ga1}), we have
\begin{eqnarray}
-\hat{\Theta}_1(t)&=&\Big(M_1'(t)-\Gamma_1'(t)\Upsilon_1^{\dag}(t)H_{1}(t)\Big)\Theta_1(t)\nonumber\\
&&
+\Big(\bar{M}_1'(t)-\Gamma_1'(t)\Upsilon_1^{\dag}(t)\bar{H}_{1}(t)\Big)\bar{\Theta}_1(t),\nonumber
\end{eqnarray}
Thus, $\Theta_1(t)$ obeys the following dynamic equation
\begin{eqnarray}
d\Theta_1(t)
&=&-\Big[\Big(M_1'(t)-\Gamma_1'(t)\Upsilon_1^{\dag}(t)H_{1}(t)\Big)\Theta_1(t)\nonumber\\
&&
+\Big(\bar{M}_1'(t)-\Gamma_1'(t)\Upsilon_1^{\dag}(t)\bar{H}_{1}(t)\Big)\bar{\Theta}_1(t)\Big]dt\nonumber\\
&&+\bar{\Theta}_1(t)dw(t)\label{c54}
\end{eqnarray}

Similar to Theorem \ref{theorem4}, we have the following results.

\begin{theorem}\label{the2}
If $Range \Big(\Gamma_1(t)\Big)\subseteq Range \Big(\Upsilon_1(t)\Big), $  then Problem (IR-LQ) is solvable if and only if
there exists $z_{1}(t)$ such that $P_1(T)x(T)=\Theta(T)=0$ where $x(t)$ obeys
\begin{eqnarray}
dx(t)&=&\Big\{\Big(M_1(t)-H_{1}(t)\Upsilon_1^{\dag}(t)\Gamma_1(t)\Big)x(t)\nonumber\\
&&+H_{1}(t)\Big(I-\Upsilon_1^{\dag}(t)\Upsilon_1(t)\Big)z_1(t)\Big\}dt
+\Big\{\Big(\bar{M}_1(t)\nonumber\\
&&-\bar{H}_{1}(t)\Upsilon_1^{\dag}(t)\Gamma_1(t)\Big)x(t)+\bar{H}_{1}(t) \Big(I-\Upsilon_1^{\dag}(t)\nonumber\\
&&\times \Upsilon_1(t)\Big)z_1(t)\Big\}dw(t).\label{c58}
\end{eqnarray} In this case,
In this case,
\begin{eqnarray}
\Theta_1(t)&=&\bar{\Theta}_1(t)=0, \nonumber \\
\Theta(t)&=&P_1(t) x(t).
\end{eqnarray}
$\bar{\Theta}(t)$ is given by (\ref{m26}),
 the controller $u(t)$ is given by (\ref{rY2}) and $u_1(t)$ is given from (\ref{m24}) as
\begin{eqnarray}
u_1(t)&=&-\Upsilon_1^{\dag}(t)\Gamma_1(t)x(t)+[I-\Upsilon_1^{\dag}(t)\Upsilon_1(t)]z_1(t),\label{YY2}
\end{eqnarray}

%
%
%
%
%
%
\end{theorem}
\emph{Proof.}  The proof is similar to Theorem \ref{theorem4}. In fact, in the case of  $Range \Big(\Gamma_1(t)\Big)\subseteq Range \Big(\Upsilon_1(t)\Big),$
if there exists $z_1(t)$ such that $P_1(T)x(T)=0,$ then $\Theta_1(T)=\Theta(T)=0.$ Considering the equation (\ref{c54}), it follows that
\begin{eqnarray}
\Theta_1(t)=0, \bar{\Theta}_1(t)=0.\nonumber
\end{eqnarray}
This implies that $\Theta(t)=P_1(t)x(t)$ and from (\ref{m26}).
Thus, 
 $u_1(t)$ given from (\ref{m24}) as in (\ref{YY2}).
Furthermore, $x(t)$ given by (\ref{m25}) is reduced to (\ref{c58}) and $\Theta(t)=P_1(t)x(t)$ solves FBSDEs (\ref{n2}), (\ref{n3}) and (\ref{c8}).
 \hfill $\blacksquare$\\
\begin{remark}
If $\Upsilon_1(t)$ is invertible, i.e., $I-\Upsilon_1^{\dag}(t)\Upsilon_1(t)=0$,  then from Theorem \ref{the2},  Problem (IR-LQ) is solvable if and only if
 $P_1(T)x(T)=\Theta(T)=0$.  In this case, the controller $u_1(t)$ is given by 
where
\begin{eqnarray}
u_1(t)&=&-\Upsilon_1^{-1}(t)\Gamma_1(t)x(t),\label{YYJ2}
\end{eqnarray}
and (\ref{c58}) is reduced to
\begin{eqnarray}
dx(t)&=&\Big\{\Big(M_1(t)-H_{1}(t)\Upsilon_1^{-1}(t)\Gamma_1(t)\Big)x(t)\Big\}dt
\nonumber \\ &&+\Big\{\Big(\bar{M}_1(t)-\bar{H}_{1}(t)\Upsilon_1^{-1}(t)\Gamma_1(t)\Big)x(t)\Big\}dw(t).\nonumber\\\label{Austr1}
\end{eqnarray}

\end{remark}



\subsection{The Second Layer}

We now consider the case of
\begin{eqnarray}
Range \Big(\Gamma_1(t)\Big)\nsubseteq Range \Big(\Upsilon_1(t)\Big),\label{m28}
\end{eqnarray} for any $P_1(T).$ It will be seen that most procedures in this subsection are similar to those in last subsection.

In view of (\ref{m28}),  it is known that $rank (\Upsilon_1(t))=m_1(t)< m-m_0(t)$, and $rank\Big(I-\Upsilon_1^{\dag}(t)\Upsilon_1(t)\Big)=m-m_0(t)-m_1(t)>0$.
Thus there exists an elementary row transformation matrix $T_1(t)$ such that
\begin{eqnarray}
T_1(t)\Big(I-\Upsilon_1^{\dag}(t)\Upsilon_1(t)\Big)=\left[
                              \begin{array}{c}
                                 0 \\
                                 \Upsilon_{T_1}(t)\\
                              \end{array}
                            \right],\label{QD1}
\end{eqnarray}
where $\Upsilon_{T_1}(t)\in R^{[m-m_0(t)-m_1(t)]\times [m-m_0(t)]}$ is full row rank. Further denote
\begin{eqnarray}
 \left[ \begin{array}{cc} \ast & C_1'(t) \\
                              \end{array} \right] &=&\Gamma'_1(t)\Big(I-\Upsilon_1^{\dag}(t)\Upsilon_1(t)\Big){T_1}^{-1}(t),\nonumber\\
\left[ \begin{array}{cc}\ast & B_1(t)\\
                              \end{array}
                            \right] &=&H_{1}(t)\Big(I-\Upsilon_1^{\dag}(t)\Upsilon_1(t)\Big){T_1}^{-1}(t),\nonumber\\
\left[ \begin{array}{cc}\ast & \bar{B}_1(t)\\
                              \end{array}
                            \right] &=& \bar{H}_{1}(t)\Big(I-\Upsilon_1^{\dag}(t)\Upsilon_1(t)\Big) {T_1}^{-1}(t).\nonumber
\end{eqnarray}

 Denote
\begin{eqnarray}
A_1(t)&=&M_1(t)-H_{1}(t)\Upsilon_1^{\dag}(t)\Gamma_{1}(t),\nonumber\\
\bar{A}_1(t)&=&\bar{M}_1(t)-\bar{H}_{1}(t)\Upsilon_1^{\dag}(t)\Gamma_{1}(t),\nonumber\\
D_1(t)&=&D_0(t)+F_0(t)\Big(I-P_1(t)\bar{F}_0(t)\Big)^{\dag} P_1(t)\nonumber\\
&&\times \bar{D}_0(t)-B_{01}(t)\Upsilon_1^{\dag}(t)B_{01}'(t),\nonumber\\
\bar{D}_1(t)&=&\bar{D}_0(t)+\bar{F}_0(t)\Big(I-P_1(t)\bar{F}_0(t)\Big)^{\dag} P_1(t)\nonumber\\
&&\times \bar{D}_0(t)-\bar{H}_{1}(t)\Upsilon_1^{\dag}(t)H_{1}'(t),\nonumber\\
F_1(t)&=&F_0(t)\Big(I-P_1(t)\bar{F}_0(t)\Big)^{\dag} -H_{1}(t)\Upsilon_1^{\dag}(t)\nonumber\\
&&\times \bar{H}_{1}'(t),\nonumber\\
\bar{F}_1(t)&=&\bar{F}_0(t)\Big(I-P_1(t)\bar{F}_0(t)\Big)^{\dag} -\bar{H}_{1}(t)\Upsilon_1^{\dag}(t)\nonumber\\
&&\times \bar{H}_{1}'(t).\nonumber
\end{eqnarray}

\begin{remark}
Based on Lemma \ref{leminv}, $D_1(t)$ and $\bar{F}_1(t)$ are symmetric, $\bar{D}_1(t)$ and $F_1(t)$ are rewritten as follows:
\begin{eqnarray}
\bar{D}_1(t)&=&\Big(I-\bar{F}_0(t)P_1(t)\Big)^{\dag} \bar{D}_0(t)-\bar{H}_{1}(t) \Upsilon_1^{\dag}(t)\nonumber\\
&&\times H_{1}'(t).\nonumber
\end{eqnarray}
\end{remark}
Similar to Theorem \ref{theorem4}, we have the following results.

\begin{theorem}\label{the3}
If $Range \Big(\Gamma_1(t)\Big)\nsubseteq Range \Big(\Upsilon_1(t)\Big)$, Problem (IR-LQ) is solvable if and only if there exists $u_2(t)\in R^{m-m_0(t)-m_1(t)}$ such that
\begin{eqnarray}
0&=&C_1(t)x(t)+B_{1}'(t)\Theta_1(t)+\bar{B}_{1}'(t)\bar{\Theta}_1(t),\label{c59}
\end{eqnarray}
where $u_2(t),x(t),\Theta_1(t)$ and $\bar{\Theta}_1(t)$ satisfy the FBSDEs:
\begin{eqnarray}
dx(t)&=&\Big[A_1(t)x(t)+D_1(t)\Theta_1(t)+F_1(t)\bar{\Theta}_1(t)\nonumber\\
&&+B_{1}(t)u_2(t)\Big]dt+\Big[\bar{A}_1(t)x(t)\nonumber\\
&&+\bar{D}_1(t)\Theta_1(t)+\bar{F}_1(t)\bar{\Theta}_1(t)+\bar{B}_{1}(t)\nonumber\\
&&\times u_{2}(t)\Big]dw(t),\label{c57}\\
d\Theta_1(t)&=&-\Big[A'_1(t)\Theta_1(t)+\bar{A}'_{1}(t)\bar{\Theta}_1(t)+C_1'(t)\nonumber\\
&&\times u_2(t)\Big]dt +\bar{\Theta}_1(t)dw(t),\label{c56}
\end{eqnarray}
with $\Theta_1(T)=-P_1(T)x(T).$

\end{theorem}

\emph{Proof.} ``Necessity"
From (\ref{c52}), one has
\begin{eqnarray}
u_1(t)&=&-\Upsilon_1^{\dag}(t)\Big(\Gamma_1(t)x(t)+H_{1}'(t)\Theta_1(t)+\bar{H}_{1}'(t)\nonumber\\
&&\times \bar{\Theta}_1(t)\Big)+\Big(I-\Upsilon_1^{\dag}(t)\Upsilon_1(t)\Big)z_1(t)\label{n4}
\end{eqnarray}
and
\begin{eqnarray}
0&=&\Big(I-\Upsilon_1(t)\Upsilon_1^{\dag}(t)\Big)\Big(\Gamma_1(t)x(t)+H_{1}'(t)\Theta_1(t)\nonumber\\
&&+\bar{H}_{1}'(t)\bar{\Theta}_1(t)\Big).\label{z8}
\end{eqnarray}
From (\ref{QD1}), we have
 \begin{eqnarray}
T_1(t)\Big(I-\Upsilon_1^{\dag}(t)\Upsilon_1(t)\Big)z_1(t)=\left[
                              \begin{array}{c}
                                 0  \\
                                 u_2(t)\\
                              \end{array}
                            \right], \label{jnj1}
\end{eqnarray}
where $u_2(t)=\Upsilon_{T_1}(t) z_1(t)$ and $\Upsilon_{T_1}(t)$ is full row rank.
 Similar to the lines as in the proof of
Theorem \ref{theorem4}, (\ref{c59}) follows directly from (\ref{z8}).
By substituting (\ref{m26}), (\ref{c48}) and (\ref{n4}) into (\ref{n3}), we have the dynamic (\ref{c57}).
From (\ref{c49}), it yields that
\begin{eqnarray}
0&=&\dot{P}_1(t)x(t)+P_1(t)\Big[A_0(t)x(t)
+D_0(t)P_1(t)\Theta_1(t)\nonumber\\
&&+F_0(t)\bar{\Theta}(t)+B_{0}(t)u_1(t)\Big]+A'_0(t)P_1(t)x(t)\nonumber\\
&&+\hat{\Theta}_1(t)+A'_0(t)\Theta_1(t)+\bar{A}'_0(t)\bar{\Theta}(t)\nonumber\\
&&+C_0'(t)u_1(t).
\end{eqnarray}
By using the similar lines on $\Theta(t)$ as in last subsection, the dynamic (\ref{c56}) of $\Theta_1$ follows.

``Sufficiency" By taking reverse procedures to the proof of Necessity, it is verified that Problem (IR LQ) is solvable.
\hfill $\blacksquare$\\

\begin{remark}
Applying (\ref{jnj1}),  we have that
\begin{eqnarray}
\Big(I-\Upsilon_1^{\dag}(t)\Upsilon_1(t)\Big)z_1(t)&=&T^{-1}_1(t)\left[
                              \begin{array}{c}
                                 0  \\
                                 u_2(t)\\
                              \end{array}
                            \right]\nonumber\\
 &=&G_1(t) u_2(t), \label{jny1}
\end{eqnarray}
where $
 \left[ \begin{array}{cc} \ast & G_1(t) \\
                              \end{array} \right] ={T_1}^{-1}(t)$.  Then  (\ref{n4}) can be rewritten as
\begin{eqnarray}
u_1(t)&=&-\Upsilon_1^{\dag}(t)\Big(\Gamma_1(t)x(t)+H_{1}'(t)\Theta_1(t)+\bar{H}_{1}'(t)\nonumber\\
&&\times \bar{\Theta}_1(t)\Big)+G_1(t)u_2(t).\label{na4}
\end{eqnarray}
\end{remark}
From Theorem \ref{the3}, the problem is now reformulated as a similar problem stated in Theorem \ref{theorem4}, which can be solved with similar lines as in Subsection A.

Now we have the similar results as in Observation \ref{obser1}.
\begin{observation}
The solution to FBSDEs (\ref{n2}), (\ref{n3}) and (\ref{c8}) is homogeneous, i.e., $\Theta_1(t) =0,$ if and only if the Riccati equation (\ref{p1}) is regular, i.e.,
$Range \Big(\Gamma_1(t)\Big)\subseteq Range \Big(\Upsilon_1(t)\Big)$, and there exists $z_1(t)$ such that $P_1(T)x(T)=0$ where the dynamic of $x(t)$
is given by (\ref{c58}).
\end{observation}
\emph{Proof.} ``Necessity" Assume that the Riccati equation (\ref{p1}) is not regular, then according to
Theorem \ref{the3}, $\Theta_1(t)$ satisfies (\ref{c56}) which is not equal to zero. This is a contradiction.

``Sufficiency" The sufficiency follows from Theorem \ref{the2}.
The proof is now completed. \hfill $\blacksquare$


Based on Theorem \ref{the3}, the solvability of Problem (IR-LQ) is converted into finding $u_2(t)$ to achieve (\ref{c59}) associated with (\ref{c57})-(\ref{c56}). The following arguments are similar to subsection A.
Let
\begin{eqnarray}
\Theta_2(t)=\Theta_1(t)-P_2(t)x(t),\label{sec1}
\end{eqnarray}
where $P_2(t)$ obeys
\begin{eqnarray}
0&=&\dot{P}_2(t)+P_2(t)A_1(t)+A_1'(t)P_2(t)+P_2(t)D_1(t)\nonumber\\
&&\times P_2(t)+\Big(\bar{A}'_1(t)+P_2(t)F_1(t)\Big)\Big(I-P_2(t)\bar{F}_{1}(t)\Big)^{\dag}\nonumber\\
&&\times P_2(t)\Big(\bar{A}_1(t)
+\bar{D}_1(t)P_2(t)\Big)-\Gamma_2'(t)\Upsilon_2^{\dag}(t)\Gamma_2(t), \nonumber\\\label{sec2}
\end{eqnarray}
where
\begin{eqnarray}
\Upsilon_2(t)&=&\bar{B}_{1}'(t)\Big(I-P_2(t)\bar{F}_{1}(t)\Big)^{\dag}P_2(t)\bar{B}_{1}(t), \label{sec3}\\
\Gamma_2(t)&=&C_1(t)+B_{1}'(t)P_2(t)+\bar{B}_{1}'(t)\Big(I-P_2(t)\nonumber\\
&&\times \bar{F}_{1}(t)\Big)^{\dag}P_2(t)\Big(\bar{A}_1(t)+\bar{D}_{1}(t) P_2(t)\Big),\label{sec4}
\end{eqnarray}
with the arbitrary terminal value $P_2(T)$.

We assume that
$$d\Theta_2(t)=\hat{\Theta}_2(t)dt+\bar{\Theta}_2(t)dw(t), $$
where $\hat{\Theta}_2(t)$ and $\bar{\Theta}_2(t)$ are to be determined. Then by taking It\^{o}'s formula to (\ref{sec1}), it is obtained that
\begin{eqnarray}
d\Theta_1(t)
&=&\dot{P}_2(t)x(t)dt+P_2(t)dx(t)+d\Theta_2(t)\nonumber\\
&=&\dot{P}_2(t)x(t)dt+P_2(t)\Big[A_1(t)x(t)
+D_1(t)\Theta_1(t)\nonumber\\
&&+F_1(t)\bar{\Theta}_1(t)+B_{1}(t)u_2(t)\Big]dt+P_2(t)\Big[\bar{A}_1(t)\nonumber\\
&&\times x(t)
+\bar{D}_1(t)\Theta_1(t)+\bar{F}_1(t)\bar{\Theta}_1(t)+\bar{B}_{1}(t)\nonumber\\
&&\times u_2(t)\Big]dw(t)+\hat{\Theta}_2(t)dt+\bar{\Theta}_2(t)dw(t).
\end{eqnarray}
Combining with (\ref{sec1}), we have
\begin{eqnarray}
0&=&\dot{P}_2(t)x(t)+P_2(t)\Big[A_1(t)x(t)
+D_1(t)\Theta_1(t)\nonumber\\
&&+F_1(t)\bar{\Theta}_1(t)+B_{1}(t)u_2(t)\Big]+A'_1(t)P_2(t)x(t)\nonumber\\
&&+\hat{\Theta}_2(t)+A'_1(t)\Theta_2(t)+\bar{A}'_1(t)\bar{\Theta}_1(t)\nonumber\\
&&+C_1'(t)u_2(t),\label{sec5}\\
\bar{\Theta}_1(t)&=&\bar{\Theta}_2(t)+P_2(t)\Big[\bar{A}_1(t)x(t)
+\bar{D}_1(t)\Theta_1(t)\nonumber\\
&&+\bar{F}_1(t)\bar{\Theta}_1(t)+\bar{B}_{1}(t)u_2(t)\Big].\label{sec6}
\end{eqnarray}
Thus it is obtained from (\ref{sec6}) that
\begin{eqnarray}
\Big(I-P_2(t)\bar{F}_{1}(t)\Big)\bar{\Theta}_1(t)
&=&\Big[\bar{\Theta}_2(t)+P_2(t)\Big(\bar{A}_1(t)x(t)\nonumber\\
&&
+\bar{D}_1(t)\Theta_1(t)+\bar{B}_{1}(t)u_2(t)\Big)\Big],\nonumber
\end{eqnarray}
which gives that
\begin{eqnarray}
\bar{\Theta}_1(t)&=&\Big(I-P_2(t)\bar{F}_{1}(t)\Big)^{\dag}\Big[\bar{\Theta}_2(t)+P_2(t)\Big(\bar{A}_1(t)x(t)\nonumber\\
&&+\bar{D}_1(t)\Theta_1(t)+\bar{B}_{1}(t)u_2(t)\Big)\Big]+\Big[I-\Big(I-P_2(t)\nonumber\\
&&\times\bar{F}_{1}(t)\Big)^{\dag}\Big(I-P_2(t)\bar{F}_{1}(t)\Big)\Big]\varphi_1(t),\label{sec7}
\end{eqnarray}
where $\varphi_1(t)$ is an arbitrary vector with compatible dimension.

Substituting (\ref{sec7}) into (\ref{c61}), it yields that
\begin{eqnarray}
0
&=&\Upsilon_0(t)u(t)
+\Gamma_0(t)x(t)+B'(t)P_1(t)x(t)+B'(t)\nonumber\\
&&\times \Theta_1(t)+\bar{B}'(t)\Big(I-P_1(t)\bar{F}_{0}(t)\Big)^{\dag}\Big(I-P_2(t)\bar{F}_{1}(t)\Big)^{\dag}\nonumber\\
&&\times \Big[\bar{\Theta}_2(t)+P_2(t)\Big(\bar{A}_1(t)x(t)+\bar{D}_1(t)\Theta_1(t)+\bar{B}_{1}(t)\nonumber\\
&&\times u_2(t)\Big)\Big]+\bar{B}'(t)\Big(I-P_1(t)\bar{F}_{0}(t)\Big)^{\dag}\Big[I-\Big(I-P_2(t)\nonumber\\
&&\times\bar{F}_{1}(t)\Big)^{\dag}\Big(I-P_2(t)\bar{F}_{1}(t)\Big)\Big]\varphi_1(t)\nonumber\\
&&+\bar{B}'(t)\Big(I-P_1(t)\bar{F}_{0}(t)\Big)^{\dag}P_1(t)\Big(\bar{A}_0(t)x(t)\nonumber\\
&&+\bar{D}_0(t)\Theta(t)+\bar{B}_{0}(t)u_1(t)\Big)\Big]
\end{eqnarray}
Similar to the derivation of (\ref{z17}), there holds that \begin{eqnarray}
0&=&\bar{B}'(t)\Big(I-P_1(t)\bar{F}_{0}(t)\Big)^{\dag}\Big[I-\Big(I-P_2(t)\bar{F}_{1}(t)\Big)^{\dag}\nonumber\\
&&\times \Big(I-P_2(t)\bar{F}_{1}(t)\Big)\Big].\nonumber\\\label{sec20}
\end{eqnarray}
Similar to Remark \ref{rem3}, we have the following results.
\begin{remark}\label{rem4}
Using (\ref{sec20}) and the definitions of $\bar{B}_0(t), \bar{D}_0(t), \bar{F}_0(t),\bar{B}_{01}(t), \bar{B}_1(t)$, it is immediate to obtain that
\begin{eqnarray}
0&=&\bar{B}_0'(t)\Big(I-P_1(t)\bar{F}_{0}(t)\Big)^{\dag}\Big[I-\Big(I-P_2(t)\bar{F}_{1}(t)\Big)^{\dag}\nonumber\\
&&\times \Big(I-P_2(t)\bar{F}_{1}(t)\Big)\Big],\nonumber\\
0&=&\bar{D}_0'(t)\Big(I-P_1(t)\bar{F}_{0}(t)\Big)^{\dag}\Big[I-\Big(I-P_2(t)\bar{F}_{1}(t)\Big)^{\dag}\nonumber\\
&&\times \Big(I-P_2(t)\bar{F}_{1}(t)\Big)\Big],\nonumber\\
0&=&\bar{F}_0(t)\Big(I-P_1(t)\bar{F}_{0}(t)\Big)^{\dag}\Big[I-\Big(I-P_2(t)\bar{F}_{1}(t)\Big)^{\dag}\nonumber\\
&&\times \Big(I-P_2(t)\bar{F}_{1}(t)\Big)\Big],\nonumber\\
0&=&\bar{B}_1'(t)\Big[I-\Big(I-P_2(t)\bar{F}_{1}(t)\Big)^{\dag}\Big(I-P_2(t)\bar{F}_{1}(t)\Big)\Big],\nonumber\\
0&=&\bar{B}_{01}'(t)\Big[I-\Big(I-P_2(t)\bar{F}_{1}(t)\Big)^{\dag}\Big(I-P_2(t)\bar{F}_{1}(t)\Big)\Big],\nonumber\\
0&=&F_1(t)\Big[I-\Big(I-P_2(t)\bar{F}_{1}(t)\Big)^{\dag}\Big(I-P_2(t)\bar{F}_{1}(t)\Big)\Big],\nonumber\\
0&=&\bar{F}_1(t)\Big[I-\Big(I-P_2(t)\bar{F}_{1}(t)\Big)^{\dag}\Big(I-P_2(t)\bar{F}_{1}(t)\Big)\Big].\nonumber
\end{eqnarray}
\end{remark}

Substituting (\ref{sec7}) into (\ref{c52}) and using the definition of $\bar{B}_{01}(t)$ and Remark 4, it is derived that
\begin{eqnarray}
0&=&\Upsilon_1(t)u_1(t)+\Gamma_1(t)x(t)+H_{1}'(t)P_2(t)x(t)+H_{1}'(t)\nonumber\\
&&\times \Theta_2(t)+\bar{H}_{1}'(t)\Big(I-P_2(t)\bar{F}_{1}(t)\Big)^{\dag}\Big[\bar{\Theta}_2(t)+P_2(t)\nonumber\\
&&\times \Big(\bar{A}_1(t)x(t)+\bar{D}_1(t)\Theta_1(t)+\bar{B}_{1}(t)u_2(t)\Big)\Big].\label{sec21}
\end{eqnarray}

Substituting (\ref{sec7}) into (\ref{c59})  yields that
\begin{eqnarray}
0&=&C_1(t)x(t)+B_{1}'(t)\Theta_1(t)+\bar{B}_{1}'(t)\bar{\Theta}_1(t)\nonumber\\
&=& \Big[\bar{B}_{1}'(t)\Big(I-P_2(t)\bar{F}_{1}(t)\Big)^{\dag} P_2(t)\bar{B}_{1}(t)\Big]u_2(t)\nonumber\\
&&+\Big[C_{1}(t)+B_{1}'(t)P_2(t)\nonumber\\
&&+\bar{B}_{1}'(t)\Big(I-P_2(t)\bar{F}_{1}(t)\Big)^{\dag}P_2(t)\Big(\bar{A}_1(t)+\bar{D}_{1}(t)\nonumber\\
&&\times P_2(t)\Big)\Big]x(t)+B_{11}'(t)\Theta_2(t)+\bar{B}_{11}'(t)\bar{\Theta}_2(t)\nonumber\\
&&
+ \bar{B}_1'(t)\Big[I-\Big(I-P_2(t)\bar{F}_{1}(t)\Big)^{\dag}\Big(I-P_2(t)\nonumber\\
&&\times\bar{F}_{1}(t)\Big)\Big]\varphi_1(t)\nonumber\\
&=&\Upsilon_2(t)u_2(t)+\Gamma_2(t)x(t)+B_{11}'(t)\Theta_2(t)\nonumber\\
&&+\bar{B}_{11}'(t)\bar{\Theta}_2(t)+\bar{B}_1'(t)\Big[I-\Big(I-P_2(t)\bar{F}_{1}(t)\Big)^{\dag}\nonumber\\
&&\times\Big(I-P_2(t)\bar{F}_{1}(t)\Big)\Big]\varphi_1(t)\nonumber\\
&=&\Upsilon_2(t)u_2(t)+\Gamma_2(t)x(t)+H_{2}'(t)\Theta_2(t)\nonumber\\
&&+\bar{H}_{2}'(t)\bar{\Theta}_2(t),\label{sec9}
\end{eqnarray}
where the matrices $\Upsilon_2(t)$ and $\Gamma_2(t)$ are respectively as (\ref{sec3}) and (\ref{sec4}), and
\begin{eqnarray}
H_{2}'(t)&=&B_{1}'(t)+\bar{B}_{1}'(t)\Big(I-P_2(t)\bar{F}_1(t)\Big)^{\dag}P_2(t)\bar{D}_1(t),\nonumber\\
\bar{H}_{2}'(t)&=&\bar{B}_{1}'(t)\Big(I-P_2(t)\bar{F}_1(t)\Big)^{\dag}.\nonumber
\end{eqnarray}
This gives that
\begin{eqnarray}
u_2(t)&=&-\Upsilon_2^{\dag}(t)\Big(\Gamma_2(t)x(t)+H_{2}'(t)\Theta_2(t)\nonumber\\
&&+\bar{H}_{2}'(t)\bar{\Theta}_2(t)\Big)+\Big(I-\Upsilon_2^{\dag}(t)\Upsilon_2(t)\Big)z_2(t),\label{sec18}
\end{eqnarray}

Similar to in Subsection A,  two cases are to be considered, one is the case of $Range \Big(\Gamma_2(t)\Big)\subseteq Range \Big(\Upsilon_2(t)\Big)$
and the other is $Range \Big(\Gamma_2(t)\Big)\nsubseteq Range \Big(\Upsilon_2(t)\Big).$

We firstly consider the case that
\begin{eqnarray}
Range \Big(\Gamma_2(t)\Big)\subseteq Range \Big(\Upsilon_2(t)\Big).\label{ga2}
\end{eqnarray}

By substituting (\ref{sec1}), (\ref{sec7}) and (\ref{sec18}) into the dynamic of the system (\ref{c57}) and using Remark \ref{rem4}, we have
\begin{eqnarray}
dx(t)&=&\Big\{\Big(M_2(t)-\hat{H}_{2}(t)\Upsilon_2^{\dag}(t)\Gamma_2(t)\Big)x(t)\nonumber\\
&&+\Big[D_1(t)+F_1(t)\Big(I-P_2(t)\bar{F}_1(t)\Big)^{\dag}P_2(t)\bar{D}_{1}(t)\nonumber\\
&&-\hat{H}_{2}(t)\Upsilon_2^{\dag}(t)H_{2}'(t)\Big]\Theta_1(t)
+\Big[F_1(t)\Big(I-P_2(t)\nonumber\\
&&\times \bar{F}_1(t)\Big)^{\dag}-\hat{H}_{2}(t)\Upsilon_2^{\dag}(t)\bar{H}_{2}'(t)\Big]\bar{\Theta}_2(t)\nonumber\\
&&+\hat{H}_{2}(t)\Big(I-\Upsilon_2^{\dag}(t)\Upsilon_2(t)\Big)z_2(t)\Big\}dt\nonumber\\
&&+\Big\{\Big(\bar{M}_2(t)-\hat{\bar{H}}_{2}(t)\Upsilon_1^{\dag}(t)\Gamma_1(t)\Big)x(t)\nonumber\\
&&+\hat{\bar{H}}_{2}(t) \Big(I-\Upsilon_2^{\dag}(t)\Upsilon_2(t)\Big)z_1(t)+\Big[\bar{D}_1(t)\nonumber\\
&&+\bar{F}_1(t)\Big(I-P_2(t)\bar{F}_1(t)\Big)^{\dag}P_2(t)\bar{D}_{1}(t)-\hat{\bar{H}}_{2}(t)\nonumber\\
&&\times\Upsilon_2^{\dag}(t)H_{2}'(t)\Big]\Theta_2(t)+\Big[\bar{F}_1(t)\Big(I-P_2(t) \bar{F}_1(t)\Big)^{\dag}\nonumber\\
&&-\hat{\bar{H}}_{2}(t)\Upsilon_2^{\dag}(t)\bar{H}_{2}'(t)\Big]
\bar{\Theta}_2(t)\Big\}dw(t),\label{sec11}
\end{eqnarray}
where
\begin{eqnarray}
M_2(t)&=&A_1(t)
+D_{1}(t)P_2(t)+F_1(t)[I-P_2(t)\bar{F}_{1}(t)]^{\dag} \nonumber\\
&&\times P_2(t)[\bar{A}_1(t)+\bar{D}_1(t)P_2(t)]\nonumber\\
\bar{M}_2(t)&=&\bar{A}_1(t)
+\bar{D}_{1}(t)P_2(t)+\bar{F}_1(t)[I-P_2(t)\bar{F}_{1}(t)]^{\dag}\nonumber\\
&&\times P_2(t)[\bar{A}_1(t)+\bar{D}_1(t)P_2(t)]\nonumber\\
\hat{H}_{2}(t)&=&B_{1}(t)+F_1(t)\Big(I-P_2(t)\bar{F}_1(t)\Big)^{\dag}P_2(t)\bar{B}_{1}(t),\nonumber\\
\hat{\bar{H}}_{2}(t)&=&\bar{B}_{1}(t)+\bar{F}_1(t)\Big(I-P_2(t)\bar{F}_1(t)\Big)^{\dag}P_2(t)\bar{B}_{1}(t).\nonumber
\end{eqnarray}
On the other hand, from (\ref{sec5}), it is further obtained that
\begin{eqnarray}
0&=&\Big[\dot{P}_2(t)+P_2(t)A_1(t)+A'_1(t)P_2(t)
+P_2(t)D_1(t)\nonumber\\
&&\times P_2(t)+\Big(\bar{A}'_1(t)+P_2(t)F_1(t)\Big)\Big(I-P_2(t)\bar{F}_{1}(t)\Big)^{\dag}\nonumber\\
&&\times P_2(t)\Big(\bar{A}_1(t)+\bar{D}_1(t)P_2(t)\Big)\Big]x(t)+\Big[A'_1(t)\nonumber\\
&&+P_2(t)D_1(t)+\Big(\bar{A}'_1(t)+P_2(t)F_1(t)\Big)\Big(I-P_2(t)\nonumber\\
&&\times \bar{F}_{1}(t)\Big)^{\dag}P_2(t)\bar{D}_1(t)\Big]\Theta_2(t)+\Big[C_1'(t)+P_2(t)B_{1}(t)\nonumber\\
&&+\Big(\bar{A}'_1(t)+P_2(t)F_1(t)\Big)\Big(I-P_2(t)\bar{F}_{1}(t)\Big)^{\dag}P_2(t)\nonumber\\
&&\times \bar{B}_{1}(t)\Big]u_2(t)+\hat{\Theta}_2(t)+\Big(\bar{A}'_1(t)+P_2(t)F_1(t)\Big)\nonumber\\
&&\times \Big(I-P_2(t)\bar{F}_{1}(t)\Big)^{\dag}\bar{\Theta}_2(t)+\Big(\bar{A}'_1(t)+P_2(t)F_1(t)\Big)\nonumber
\end{eqnarray}
\begin{eqnarray}
&&\times\Big[I-\Big(I-P_2(t)\bar{F}_{1}(t)\Big)^{\dag}\Big(I-P_2(t)\bar{F}_{1}(t)\Big)\Big]\varphi_1(t)\nonumber
\end{eqnarray}
By denoting
\begin{eqnarray}
\hat{M}_2(t)&=&A_1(t)
+D_{1}(t)P_2(t)+F_1(t)P_2(t)[I-\bar{F}_{1}(t)\nonumber\\
&&\times P_2(t)]^{\dag}  [\bar{A}_1(t)+\bar{D}_1(t)P_2(t)],\nonumber\\
\hat{\bar{M}}_2(t)&=&[I-\bar{F}_{1}(t)P_2(t)]^{\dag}[\bar{A}_1(t)+\bar{D}_1(t)P_2(t)],\nonumber\\
\hat{\Gamma}_2(t)&=&C_1(t)+B_{1}'(t)P_2(t)+\bar{B}_{1}'(t)P_2(t)\Big(I-P_2(t)\nonumber\\
&&\times \bar{F}_{1}(t)\Big)^{\dag}\Big(\bar{A}_1(t)+\bar{D}_{1}(t) P_2(t)\Big).\nonumber
\end{eqnarray}
Combining with (\ref{sec18}), the above equation becomes
\begin{eqnarray}
0&=&\Big[\dot{P}_2(t)+P_2(t)A_1(t)+A'_1(t)P_2(t)
+P_2(t)D_1(t)\nonumber\\
&&\times P_2(t)+\Big(\bar{A}'_1(t)+P_2(t)F_1(t)\Big)\Big(I-P_2(t)\bar{F}_{1}(t)\Big)^{\dag}\nonumber\\
&&\times P_2(t)\Big(\bar{A}_1(t)+\bar{D}_1(t)P_2(t)\Big)-\hat{\Gamma}_2'(t)\Upsilon_2^{\dag}(t)\Gamma_2(t)\Big]\nonumber\\
&&\times x(t)+\Big(\hat{M}_2'(t)-\hat{\Gamma}_2'(t)\Upsilon_2^{\dag}(t)H_2'(t)\Big)\Theta_2(t)\nonumber\\
&&+\hat{\Gamma}_2'(t)u_2(t)+\hat{\Theta}_2(t)+\Big(\hat{\bar{M}}_2'(t)-\hat{\Gamma}_2'(t)\Upsilon_2^{\dag}(t)\nonumber\\
&&\times \bar{H}_2'(t)\Big)\bar{\Theta}_2(t)+\hat{\Gamma}_2'(t)\Big(I-\Upsilon_2^{\dag}(t)\Upsilon_2(t)\Big)z_2(t)\nonumber\\
&&+\bar{A}'_1(t)\Big[I-\Big(I-P_2(t)\bar{F}_{1}(t)\Big)^{\dag}\nonumber\\
&&\times \Big(I-P_2(t)\bar{F}_{1}(t)\Big)\Big]\varphi_1(t), \label{sec15}
\end{eqnarray}
where Remark \ref{rem4} has been used in the last equality.
Similar to (\ref{z17}), it follows that
\begin{eqnarray}
0&=&\bar{A}'_1(t)\Big[I-\Big(I-P_2(t)\bar{F}_{1}(t)\Big)^{\dag}\Big(I-P_2(t)\bar{F}_{1}(t)\Big)\Big].\nonumber\\\label{sec16}
\end{eqnarray}
The following Lemma is similar to Lemma \ref{leminv}.
\begin{lemma}\label{leminvsec1}
Under (\ref{sec20}) and (\ref{sec16}), it holds that
\begin{enumerate}
  \item
  \begin{eqnarray}
  &&L_{1}'(t)P_2(t)\Big(I-\bar{F}_{1}(t)P_2(t)\Big)^{\dag}N_{1}(t)\nonumber\\
  &=&L_{1}'(t)\Big(I-P_2(t)\bar{F}_{1}(t)\Big)^{\dag}P_2(t)N_{1}(t),\nonumber
  \end{eqnarray}
  where $L_{1}(t),N_{1}(t)=\bar{B}_{1}(t),\bar{D}_{1}(t),\bar{A}_{1}(t),\bar{F}_{1}(t).$
  \item \begin{eqnarray}
  &&\Big(I-\bar{F}_{1}(t)P_2(t)\Big)^{\dag}L_{1}(t)\nonumber\\
  &=&\Big[I+\bar{F}_{1}(t)\Big(I-P_2(t)\bar{F}_{1}(t)\Big)^{\dag}P_2(t)\Big]L_{1}(t),\nonumber
  \end{eqnarray} where $L_{1}(t)=\bar{B}_{1}(t),\bar{D}_{1}(t),\bar{A}_{1}(t),\bar{F}_{1}(t).$
\end{enumerate}
\end{lemma}
\emph{Proof.} Based on Remark \ref{rem4}, there holds that
\begin{eqnarray}
0&=&L_{1}'(t)\Big[I-\Big(I-P_2(t)\bar{F}_{1}(t)\Big)^{\dag}\nonumber\\
&&\times\Big(I-P_2(t)\bar{F}_{1}(t)\Big)\Big],\nonumber
\end{eqnarray}
where $L_{1}(t)=\bar{A}_{1}(t), \bar{B}_{1}'(t),\bar{D}_{1}'(t),\bar{F}_{1}(t).$ Combining with the discussions in Lemma \ref{leminv},
the results follow. \hfill $\blacksquare$

\begin{lemma}\label{leminvsec2}
Under (\ref{sec20}) and (\ref{sec16}), it holds for $L_0(t)=\bar{B}_0(t),\bar{A}_0(t),\bar{D}_0(t),\bar{F}_0(t)$ that
\begin{enumerate}
  \item \begin{eqnarray}
L_0'(t)&=&L_0'(t)\Big(I-P_1(t)\bar{F}_0(t)\Big)^{\dag}\Big[I-P_2(t)\bar{F}_0(t)\nonumber\\
&&\times\Big(I-P_1(t)\bar{F}_0(t)\Big)^{\dag}\Big]^{\dag}
\Big[I-P_1(t)\bar{F}_0(t)\nonumber\\
&&-P_2(t)\bar{F}_0(t)\Big],\label{z14}
\end{eqnarray}
  \item \begin{eqnarray}
&&\Big(P_1(t)+P_2(t)\Big)L_0(t)\nonumber\\
&=&\Big[I-P_1(t)\bar{F}_0(t)-P_2(t)\bar{F}_0(t)\Big]\Big[I-P_1(t)\bar{F}_0(t)\nonumber\\
&&-P_2(t)\bar{F}_0(t)\Big]^{\dag}
\Big(P_1(t)+P_2(t)\Big)L_0(t).\label{z21}
\end{eqnarray}
\end{enumerate}
\end{lemma}
\emph{Proof.} We only state the proof for the case of $L_0(t)=\bar{B}_0(t).$
The proof for the other cases are the same.
\begin{enumerate}
  \item Using Remark \ref{rem4}, it is obtained that
\begin{eqnarray}
&&\hspace{-3mm}\bar{B}_{0}'(t)\Big(I-P_1(t)\bar{F}_0(t)\Big)^{\dag}\Big[I-P_2(t)\bar{F}_0(t)\Big(I-P_1(t)\nonumber\\
&&\hspace{-3mm}\times\bar{F}_0(t)\Big)^{\dag}\Big]^{\dag}
\Big[I-P_2(t)\bar{F}_0(t)\Big(I-P_1(t)\bar{F}_0(t)\Big)^{\dag}\Big]\nonumber\\
&=&\bar{B}_{0}'(t)\Big(I-P_1(t)\bar{F}_0(t)\Big)^{\dag}.\label{z22}
\end{eqnarray}
Multiplying $I-P_1(t)\bar{F}_0(t)$ from the right side yields that
\begin{eqnarray}
&&\bar{B}_{0}'(t)\Big(I-P_1(t)\bar{F}_0(t)\Big)^{\dag}\Big[I-P_2(t)\bar{F}_0(t)\Big(I\nonumber\\
&&-P_1(t) \bar{F}_0(t)\Big)^{\dag}\Big]^{\dag}
\Big[I-P_2(t)\bar{F}_0(t)\Big(I-P_1(t)\nonumber\\
&&\times\bar{F}_0(t)\Big)^{\dag}\Big]\Big(I-P_1(t)\bar{F}_0(t)\Big)\nonumber\\
&=&\bar{B}_{0}'(t)\Big(I-P_1(t)\bar{F}_0(t)\Big)^{\dag}\Big(I-P_1(t)\bar{F}_0(t)\Big).\nonumber
\end{eqnarray}
By using (\ref{z1}), (\ref{z14}) follows.

  \item From (\ref{z14}), it yields that
  \begin{eqnarray}
L_0(t)&=&\Big[I-\bar{F}_0(t)P_1(t)-\bar{F}_0(t)P_2(t)\Big]
\Big[I\nonumber\\
&&-\Big(I-\bar{F}_0(t)P_1(t)\Big)^{\dag}\bar{F}_0(t)P_2(t)\Big]^{\dag}\nonumber\\
&&\times\Big(I-\bar{F}_0(t)P_1(t)\Big)^{\dag}L_0(t),
\label{z11}
\end{eqnarray}
which can be written as
\begin{eqnarray}
\Big[I-\bar{F}_0(t)P_1(t)-\bar{F}_0(t)P_2(t)\Big]y_1(t)=L_{0}(t), \nonumber
\end{eqnarray}
where
\begin{eqnarray}
y_1(t)&=&\Big[I-\Big(I-\bar{F}_0(t)P_1(t)\Big)^{\dag}\bar{F}_0(t)P_2(t)\Big]^{\dag}\nonumber\\
&&\times\Big(I-\bar{F}_0(t)P_1(t)\Big)^{\dag}L_{0}(t),\nonumber
\end{eqnarray}
is a solution to the above equation. Thus we have
\begin{eqnarray}
&&Range\Big(L_{0}(t)\Big)\nonumber\\
&\subseteq& Range \Big(I-\bar{F}_0(t)P_1(t)-\bar{F}_0(t)P_2(t)\Big).\nonumber
\end{eqnarray}
This gives that
\begin{eqnarray}
L_{0}(t)&=&\Big[I-\bar{F}_0(t)P_1(t)-\bar{F}_0(t)P_2(t)\Big]\Big[I-\bar{F}_0(t)\nonumber\\
&&\times P_1(t)-\bar{F}_0(t)P_2(t)\Big]^{\dag}L_{0}(t).\nonumber
\end{eqnarray}
By multiplying $P_1(t)+P_2(t)$ from the left side to the above equation, it yields that
\begin{eqnarray}
&&\Big(P_1(t)+P_2(t)\Big)L_{0}(t)\nonumber\\
&=&\Big[I-P_1(t)\bar{F}_0(t)-P_2(t)\bar{F}_0(t)\Big]\Big(P_1(t)+P_2(t)\Big)\nonumber\\
&&\times\Big[I-\bar{F}_0(t)P_1(t)-\bar{F}_0(t)P_2(t)\Big]^{\dag}L_{0}(t).\label{z12}
\end{eqnarray}
Denote
\begin{eqnarray}
y_2(t)&=&\Big(P_1(t)+P_2(t)\Big)\Big[I-\bar{F}_0(t)P_1(t)\nonumber\\
&&-\bar{F}_0(t)P_2(t)\Big]^{\dag}L_{0}(t),\nonumber
\end{eqnarray}
then (\ref{z12}) is reduced to
\begin{eqnarray}
&&\Big[I-P_1(t)\bar{F}_0(t)-P_2(t)\bar{F}_0(t)\Big]y_2(t)\nonumber\\
&=&\Big(P_1(t)+P_2(t)\Big)L_{0}(t),\nonumber
\end{eqnarray}
thus
\begin{eqnarray}
&&Range\Big[\Big(P_1(t)+P_2(t)\Big)L_{0}(t)\Big]\nonumber\\
&\subseteq& Range \Big(I-P_1(t)\bar{F}_0(t)-P_2(t)\bar{F}_0(t)\Big).\nonumber
\end{eqnarray}
This implies that (\ref{z21}) holds. The proof is now completed. \hfill $\blacksquare$

\end{enumerate}

\begin{remark}\label{rem2}
Based on Lemma \ref{leminvsec1}, we have that
\begin{eqnarray}
\hat{H}_{2}(t)
&=&H_{2}(t)\nonumber\\
\hat{\bar{H}}_{2}(t)
&=&\bar{H}_{2}(t),\nonumber\\
\hat{M}_2(t)&=&M_2(t)\nonumber\\
\hat{\bar{M}}_2(t)&=&\bar{M}_2(t),\nonumber\\
\hat{\Gamma}_2(t)&=&\Gamma_2(t),\nonumber\\
\bar{H}_{2}'(t)
&=&\bar{B}_{1}'(t)+\bar{B}_{1}'(t)\Big(I-P_2(t)\bar{F}_{1}(t)\Big)^{\dag}P_2(t)\bar{F}_1(t),\nonumber\\
\bar{M}_2'(t)&=&[\bar{A}_1(t)+\bar{D}_1(t)P_2(t)]'[I-P_2(t)\bar{F}_{1}(t)]^{\dag}.\nonumber
\end{eqnarray}

\end{remark}

With the similar discussions on (\ref{c60}) in Subsection A, by combining with (\ref{sec2}), (\ref{sec18}), (\ref{ga2}) and Remark \ref{rem2}, we know that (\ref{sec15}) is reduced to
\begin{eqnarray}
0&=&\Big(M_2'(t)-\Gamma_2'(t)\Upsilon_2^{\dag}(t)H_{2}(t)\Big)\Theta_2(t)\nonumber\\
&&
+\Big(\bar{M}_2'(t)-\Gamma_2'(t)\Upsilon_2^{\dag}(t)\bar{H}_{2}(t)\Big)\bar{\Theta}_2(t)\nonumber\\
&&+\hat{\Theta}_2(t),\nonumber
\end{eqnarray}
that is, $\Theta_2(t)$ obeys the following dynamic equation:
\begin{eqnarray}
d\Theta_2(t)&=&-\Big[\Big(M_2'(t)-\Gamma_2'(t)\Upsilon_2^{\dag}(t)H_{2}(t)\Big)\Theta_2(t)\nonumber\\
&&
+\Big(\bar{M}_2'(t)-\Gamma_2'(t)\Upsilon_2^{\dag}(t)\bar{H}_{2}(t)\Big)\bar{\Theta}_2(t)\Big]dt\nonumber\\
&&+\bar{\Theta}_2(t)dw(t).\label{sec19}
\end{eqnarray}

The following Theorem is similar to Theorem \ref{the2}.
\begin{theorem}\label{the4}
If $Range \Big(\Gamma_2(t)\Big)\subseteq Range \Big(\Upsilon_2(t)\Big), $  then Problem (IR-LQ) is solvable if and only if
there exists $z_{2}(t)$ such that $\Big(P_1(T)+P_2(T)\Big)x(T)=0$ where $x(t)$ obeys
\begin{eqnarray}
dx(t)&=&\Big\{\Big(M_2(t)-H_{2}(t)\Upsilon_2^{\dag}(t)\Gamma_2(t)\Big)x(t)+H_{2}(t)\nonumber\\
&&\times \Big(I-\Upsilon_2^{\dag}(t)\Upsilon_2(t)\Big)z_2(t)\Big\}dt+\Big\{\Big(\bar{M}_2(t)\nonumber\\
&&-\bar{H}_{2}(t)\Upsilon_2^{\dag}(t)\Gamma_2(t)\Big)x(t)+\bar{H}_{2}(t) \Big(I-\Upsilon_2^{\dag}(t)\nonumber\\
&&\times \Upsilon_2(t)\Big)z_2(t)\Big\}dw(t).
\end{eqnarray}

In this case,
\begin{eqnarray}
\Theta_2(t)&=&\bar{\Theta}_2(t)=0, \nonumber \\
\Theta_1(t)&=&P_2(t) x(t),\nonumber \\
\Theta(t)&=&\Big(P_1(t)+P_2(t)\Big) x(t).
\end{eqnarray}
$\bar{\Theta}(t)$ and  $\bar{\Theta}_1(t)$ are respectively given by (\ref{m26}) and (\ref{sec7}). The controller $u(t)$ and $u_1(t)$ are respectively given by (\ref{rY2}) and (\ref{na4}),  and $u_2(t)$ is as
\begin{eqnarray}
u_2(t)&=&-\Upsilon_2^{\dag}(t)\Gamma_2(t)x(t)+\Big(I-\Upsilon_2^{\dag}(t)\Upsilon_2(t)\Big)z_2(t).\nonumber\\\label{sec22}
\end{eqnarray}
\end{theorem}
\emph{Proof.}  The proof is similar to Theorem \ref{the2} and thus omitted. 

\hfill $\blacksquare$\\


\subsection{The $k$th Layer}
In this subsection, we consider the general case when $$Range \Big(\Gamma_{k-1}(t)\Big)\not\subseteq Range \Big(\Upsilon_{k-1}(t)\Big). $$

Recalling the first layer in Subsection A and the second layer in Subsection B, we know that the $kth$ layer optimization can be reduced to
seek $u_{k}(t)$  satisfy the following equilibrium condition
\begin{eqnarray}
0&=&C_{k-1}x(t)+B_{k-1}'(t)\Theta_{k-1}(t)\nonumber\\
&&+\bar{B}_{k-1}'(t)\bar{\Theta}_{k-1}(t)\label{c41}
\end{eqnarray}
where $\Theta_{k-1}(t)$ and $x(t)$ are given by
\begin{eqnarray}
dx(t)
&=& \Big[A_{k-1}(t)x(t)+D_{k-1}(t)\Theta_{k-1}(t)\nonumber\\
&&+F_{k-1}(t)\bar{\Theta}_{k-1}(t)+B_{k-1}(t)\nonumber\\
&&\times u_{k}(t)\Big]dt+\Big[\bar{A}_{k-1}(t)x(t)\nonumber\\
&&+\bar{D}_k(t)\Theta_{k-1}(t)+\bar{F}_{k-1}(t)\nonumber\\
&&\times\bar{\Theta}_{k-1}(t)+\bar{B}_{k-1}(t)u_{k}(t)\Big]dw(t)\nonumber\\\label{c36}\\
d\Theta_{k-1}(t)&=&-[A_{k-1}'(t)\Theta_{k-1}(t)+\bar{A}_{k-1}'(t)\bar{\Theta}_{k-1}(t)\nonumber\\
&&+C_{k-1}'(t)u_{k}(t)]dt+\bar{\Theta}_{k-1}(t)dw(t),  \nonumber\\\label{c37}
\end{eqnarray}
where $A_{k-1}, \bar{A}_{k-1}, B_{k-1}, \bar{B}_{k-1}, C_{k-1}, D_{k-1}, \bar{D}_{k-1}, F_{k-1},\\ \bar{F}_{k-1}$ have the similar definitions as in Subsection B.  Similar to Subsection B,
it is known that $rank (\Upsilon_{k-1}(t))=m_{k-1}(t)< m-\sum_{j=0}^{k-2}m_j(t)=dim\Big(u_{k-1}\Big)$, and $rank\Big(I-\Upsilon_{k-1}^{\dag}(t)\Upsilon_{k-1}(t)\Big)=m-\sum_{j=0}^{k-1}m_j(t)>0$.
Suppose  there exists an elementary row transformation matrix $T_{k-1}(t)$ such that
\begin{eqnarray}
T_{k-1}(t)\Big(I-\Upsilon_{k-1}^{\dag}(t)\Upsilon_{k-1}(t)\Big)=\left[
                              \begin{array}{c}
                                 0 \\
                                 \Upsilon_{T_{k-1}}(t)\\
                              \end{array}
                            \right],\label{QD1}
\end{eqnarray}
where $\Upsilon_{T_{k-1}}(t)\in R^{[m-\sum_{j=0}^{k-1} m_j(t)]\times [m-\sum_{j=0}^{k-2} m_j(t)]}$ is full row rank. Further denote
\begin{eqnarray}
 \left[ \begin{array}{cc} \ast & C_{k-1}'(t) \\
                              \end{array} \right] &=&\Gamma'_{k-1}(t)\Big(I-\Upsilon_{k-1}^{\dag}(t)\Upsilon_{k-1}(t)\Big)\nonumber\\
                              &&\times {T^{-1}_{k-1}}(t),\nonumber\\
\left[ \begin{array}{cc}\ast & B_{k-1}(t)\\
                              \end{array}
                            \right] &=&H_{k-1}(t)\Big(I-\Upsilon_{k-1}^{\dag}(t)\Upsilon_{k-1}(t)\Big)\nonumber\\
                            &&\times {T^{-1}_{k-1}}(t),\nonumber\\
\left[ \begin{array}{cc}\ast & \bar{B}_{k-1}(t)\\
                              \end{array}
                            \right] &=& \bar{H}_{k-1}(t)\Big(I-\Upsilon_{k-1}^{\dag}(t)\Upsilon_{k-1}(t)\Big)\nonumber\\
                            &&\times  T_{k-1}^{-1}(t),\nonumber\\
 \left[ \begin{array}{cc} \ast & G_{k-1}(t) \\
                              \end{array} \right] &=& T^{-1}_{k-1}(t),
\end{eqnarray}
where
\begin{eqnarray}
H_{k}'(t)&=&B_{k-1}'(t)+\bar{B}_{k-1}'(t)\Big(I-P_k(t)\bar{F}_{k-1}(t)\Big)^{\dag}\nonumber\\
&&\times P_k(t)\bar{D}_{k-1}(t),\nonumber\\
\bar{H}_{k}'(t)&=&\bar{B}_{k-1}'(t)\Big(I-P_k(t)\bar{F}_{k-1}(t)\Big)^{\dag}.\nonumber
\end{eqnarray}

The associated Riccati equation is as
\begin{eqnarray}
0&=&\dot{P}_k(t)+P_k(t)A_{k-1}(t)+P_k(t)D_{k-1}(t)P_k(t)\nonumber\\
&&+A_{k-1}'(t)P_k(t)+\Big(\bar{A}_{k-1}'(t)+P_k(t)F_{k-1}(t)\Big)\nonumber\\
&&\times\Big(I-P_k(t)\bar{F}_{k-1}(t)\Big)^{\dag} P_k(t)\Big(\bar{A}_{k-1}(t)\nonumber\\
&&+\bar{D}_{k-1}(t)P_k(t)\Big)-\Gamma_k'(t)\Upsilon_k^{\dag}(t)\Gamma_k(t), \nonumber\\\label{c45}
\end{eqnarray}
with any value $P_k(T)$ and
\begin{eqnarray}
\Upsilon_k(t)&=&\bar{B}_{k-1}'(t)\Big(I-P_k(t)\bar{F}_{k-1}(t)\Big)^{\dag} P_k(t)\bar{B}_{k-1}(t)\nonumber\\
\Gamma_k(t)&=&C_{k-1}(t)+B_{k-1}'(t)P_k(t)+\bar{B}_{k-1}'(t)\Big(I-P_k(t)\nonumber\\
&&\times\bar{F}_{k-1}(t)\Big)^{\dag}P_k(t)\Big(\bar{A}_{k-1}(t)+\bar{D}_{k-1}(t)P_k(t)\Big). \nonumber
\end{eqnarray}

Similarly, we denote
\begin{eqnarray}
M_{k}(t)&=&A_{k-1}(t)+D_{k-1}'(t)P_{k}(t)+F_{k-1}(t)[I\nonumber\\
&&-P_{k}(t)\bar{F}_{k-1}(t)]^{\dag}P_k(t)[\bar{A}_{k-1}(t)\nonumber\\
&&+\bar{D}_{k-1}(t)P_{k}(t)]\label{z19}\\
\bar{M}_{k}(t)&=&\bar{A}_{k-1}(t)
+\bar{D}_{k-1}(t)P_k(t)+\bar{F}_{k-1}(t)[I-P_k(t)\nonumber\\
&&\times\bar{F}_{k-1}(t)]^{\dag} P_k(t)[\bar{A}_{k-1}(t)+\bar{D}_{k-1}(t)P_k(t)]. \nonumber\\\label{z20}
\end{eqnarray}

The following conditions are necessary for the solvability of the Problem (IR-LQ)
\begin{eqnarray}
0&=&L_{k-2}'(t)\Big(I-P_{k-1}(t)\bar{F}_{k-2}(t)\Big)^{\dag}\Big[I-\Big(I-P_k(t)\nonumber\\
&&\times \bar{F}_{k-1}(t)\Big)^{\dag}\Big(I-P_k(t)\bar{F}_{k-1}(t)\Big)\Big],\label{JNJ5}
\end{eqnarray}
and
\begin{eqnarray}
0&=&A_{k-1}'(t)\Big[I-\Big(I-P_k(t)\bar{F}_{k-1}(t)\Big)^{\dag}\nonumber\\
&&\times \Big(I-P_k(t)\bar{F}_{k-1}(t)\Big)\Big]
\end{eqnarray}
where $L_{k-2}(t)=\bar{B}_{k-2}'(t),\bar{D}_{k-2}'(t),\bar{F}_{k-2}(t).$

Taking similar discussions to Theorem \ref{the2}, we have the following result.
\begin{theorem}\label{theorem5}
If $Range \Big(\Gamma_k(t)\Big)\subseteq Range \Big(\Upsilon_k(t)\Big),$ then Problem (IR-LQ) is solvable if and only if
there exists $z_{k}(t)$ such that $\sum_{i=1}^kP_i(T)x(T)=0$ where $x(t)$ obeys the following dynamic

\begin{eqnarray}
dx(t)&=&\Big[\Big(M_k(t)-H_{k}(t)\Upsilon_k^{\dag}(t)\Gamma_k(t)\Big)x(t)\nonumber\\
&&+H_{k}(t)\Big(I-\Upsilon_k^{\dag}(t)\Upsilon_k(t)\Big)z_{k}(t)\Big]dt\nonumber\\
&&+\Big[\Big(\bar{M}_k(t)-\bar{H}_{k}(t)\Upsilon_k^{\dag}(t)\Gamma_k(t)\Big)x(t)\nonumber\\
&&+\bar{H}_{k}(t)\Big(I-\Upsilon_k^{\dag}(t) \Upsilon_k(t)\Big)z_{k}(t)\Big]dw(t),\nonumber
\end{eqnarray}
where $M_k(t)$ and $\bar{M}_k(t)$ are as defined in (\ref{z19}) and (\ref{z20}). In this case,
\begin{eqnarray}
\Theta_k(t)&=&\bar{\Theta}_k(t)=0, \nonumber \\
\Theta_i(t)&=&\Big[\sum_{j=i+1}^kP_j(t)\Big] x(t), ~~0\leq i\leq k-1,
\end{eqnarray}
where $\bar{\Theta}_i(t)$ ($0\leq i\leq k-1$) is given as
 \begin{eqnarray}
\bar{\Theta}_{i}(t)&=&\Big(I-P_{i+1}(t)\bar{F}_{i}(t)\Big)^{\dag}\Big[\bar{\Theta}_{i+1}(t)+P_{i+1}(t)\Big(\bar{A}_{i}(t)\nonumber\\
&&\times x(t)+\bar{D}_{i}(t)\Theta_{i}(t)+\bar{B}_{i}(t)u_{i+1}(t)\Big)\Big]\nonumber\\
&&+\Big[I-\Big(I-P_{i+1}(t)\bar{F}_{i}(t)\Big)^{\dag}\Big(I-P_{i+1}(t)\bar{F}_{i}(t)\Big)\Big]\nonumber\\
&&\times\varphi_{i}(t).\label{seci}
\end{eqnarray}
The controllers are given as
\begin{eqnarray}
u_i(t)&=&-\Upsilon_i^{\dag}(t)\Big(\Gamma_i(t)x(t)+H_{i}'(t)\Theta_i(t)+\bar{H}_{i}'(t)\nonumber\\
&&\times \bar{\Theta}_i(t)\Big)+G_i(t)u_{i+1}(t), ~~0\leq i\leq k-1, \nonumber\\
u_{k}(t)&=&-\Upsilon_k^{\dag}(t)\Gamma_k(t)x(t)+[I-\Upsilon_k^{\dag}(t)\Upsilon_k(t)]z_{k}(t).\nonumber
\end{eqnarray}
In the above, $\Theta_0(t)=\Theta(t)$, $\bar{\Theta}_0(t)=\bar{\Theta}(t)$, $H_0(t)=B(t)$, $\bar{H}_0(t)=\bar{B}(t)$, $u_0(t)=u(t)$.


\end{theorem}

%

Now we present one of the key results for multi-layer optimization for the Problem (IR-LQ).

\begin{theorem}\label{JNB1}
If there exists $k\geq 1$ that $\Upsilon_k(t)=0$ and $\Gamma_k(t)\neq0$ for any $P_i(T),$ ($ 1 \leq i \leq k$) then Problem (IR-LQ) is unsolvable.
\end{theorem}

\emph{Proof.} Without losses of generality, we consider the case of $k=1$.  If $\Upsilon_1(t)=0,$ $\Gamma_1(t)\not=0$,  then (\ref{p1}) becomes
\begin{eqnarray}
0&=&\dot{P}_1(t)+P_1(t)A_0(t)+A_0'(t)P_1(t)+P_1(t)D_0(t)\nonumber\\
&&\times P_1(t)+\Big(\bar{A}'_0(t)+P_1(t)F_0(t)\Big)\Big(I-P_1(t)\bar{F}_{0}(t)\Big)^{\dag}\nonumber\\
&&\times P_1(t)\Big(\bar{A}_0(t)
+\bar{D}_0(t)P_1(t)\Big).\label{z16}
\end{eqnarray}

In what follows, we will show $\Upsilon_2(t)=0$ and $\Gamma_2(t)\not=0$, where $\Upsilon_2(t)$ and $\Gamma_2(t)$ are as defined in  (\ref{sec2}).  Following (\ref{sec1}),
we have $$\Theta_1(t)=P_2(t)x(t)+\Theta_2(t),$$  where $P_2(t)$ obeys Riccati equation (\ref{sec2}).

%

In view of (\ref{z22}), one has
\begin{eqnarray}
&&\bar{B}_{0}'(t)\Big(I-P_1(t)\bar{F}_0(t)\Big)^{\dag}\Big[I-P_2(t)\bar{F}_0(t)\Big(I-P_1(t)\nonumber\\
&&\times \bar{F}_0(t)\Big)^{\dag}\Big]^{\dag}
P_2(t)\bar{F}_0(t)\Big(I-P_1(t)\bar{F}_0(t)\Big)^{\dag}\nonumber\\
&=&\bar{B}_{0}'(t)\Big(I-P_1(t)\bar{F}_0(t)\Big)^{\dag}\Big[I-P_2(t)\bar{F}_0(t)\Big(I-P_1(t)\nonumber\\
&&\times\bar{F}_0(t)\Big)^{\dag}\Big]^{\dag}
-\bar{B}_{0}'(t)\Big(I-P_1(t)\bar{F}_0(t)\Big)^{\dag}.\nonumber
\end{eqnarray}
Together with $2)$ in Lemma \ref{leminv}, $\Upsilon_2(t)$ can be reformulated as
\begin{eqnarray}
&&\Upsilon_2(t)\nonumber\\
&=&\bar{B}_{0}'(t)\Big(I-P_1(t)\bar{F}_0(t)\Big)^{\dag}\Big[I-P_2(t)\bar{F}_0(t)\Big(I-P_1(t)\nonumber\\
&&\times \bar{F}_0(t)\Big)^{\dag}\Big]^{\dag} P_2(t)
\Big[\bar{B}_{0}(t)+\bar{F}_0(t)\Big(I-P_1(t)\bar{F}_0(t)\Big)^{\dag}\nonumber\\
&&\times P_1(t)\bar{B}_{0}(t)\Big]\nonumber\\
&=&\bar{B}_{0}'(t)\Big(I-P_1(t)\bar{F}_0(t)\Big)^{\dag}\Big[I-P_2(t)\bar{F}_0(t)\Big(I-P_1(t)\nonumber\\
&&\times\bar{F}_0(t)\Big)^{\dag}\Big]^{\dag} P_2(t)\bar{B}_{0}(t)+\Big\{\bar{B}_{0}'(t)\Big(I-P_1(t)\bar{F}_0(t)\Big)^{\dag}\nonumber\\
&&\times \Big[I-P_2(t)\bar{F}_0(t)\Big(I-P_1(t)\bar{F}_0(t)\Big)^{\dag}\Big]^{\dag}\nonumber\\
&&-\bar{B}_{0}'(t)\Big(I-P_1(t)\bar{F}_0(t)\Big)^{\dag}\Big\}P_1(t)\bar{B}_{0}(t)\nonumber\\
&=&\bar{B}_{0}'(t)\Big(I-P_1(t)\bar{F}_0(t)\Big)^{\dag}\Big[I-P_2(t)\bar{F}_0(t)\Big(I-P_1(t)\nonumber\\
&&\times \bar{F}_0(t)\Big)^{\dag}\Big]^{\dag}\Big(P_1(t)+ P_2(t)\Big)\bar{B}_{0}(t)\nonumber\\
&&
-\bar{B}_{0}'(t)\Big(I-P_1(t)\bar{F}_0(t)\Big)^{\dag}P_1(t)\bar{B}_{0}(t)\nonumber\\
&=&\bar{B}_{0}'(t)\Big(I-P_1(t)\bar{F}_{0}(t)\Big)^{\dag}\Big[I-P_2(t)\bar{F}_{0}(t)\Big(I-P_1(t)\nonumber\\
&&\times\bar{F}_{0}(t)\Big)^{\dag}\Big]^{\dag} \Big(P_1(t)+P_2(t)\Big)\bar{B}_{0}(t)\label{z13}
\end{eqnarray}
where $\Upsilon_1(t)=\bar{B}_{0}'(t)\Big(I-P_1(t)\bar{F}_0(t)\Big)^{\dag}P_1(t)\bar{B}_{0}(t)=0$ has been used in the derivation of the last equality.
Combining with (\ref{z14}) in Lemma \ref{leminvsec2}, (\ref{z13}) is further rewritten as
\begin{eqnarray}
&&\Upsilon_2(t)\nonumber\\
&=&\bar{B}_{0}'(t)\Big(I-P_1(t)\bar{F}_{0}(t)\Big)^{\dag}\Big[I-P_2(t)\bar{F}_{0}(t)\Big(I-P_1(t)\nonumber\\
&&\times \bar{F}_{0}(t)\Big)^{\dag}\Big]^{\dag} \Big[I-P_1(t)\bar{F}_0(t)-P_2(t)\bar{F}_0(t)\Big]\Big[I\nonumber\\
&&-P_1(t)\bar{F}_0(t)-P_2(t)\bar{F}_0(t)\Big]^{\dag}
\Big(P_1(t)+P_2(t)\Big)\bar{B}_{0}(t)\nonumber\\
&=&\bar{B}_{0}'(t)\Big[I-{\cal P}(t) \bar{F}_0(t)\Big]^{\dag}
{\cal P}(t)\bar{B}_{0}(t).
\end{eqnarray}
where ${\cal P}(t)=P_1(t)+P_2(t)$.
Similarly, we can show that
\begin{eqnarray}
\Gamma_2(t)&=&C_0(t)+B_{0}'(t){\cal P}(t)+\bar{B}_{0}'(t)\Big(I-{\cal P}(t)\bar{F}_{0}(t)\Big)^{\dag}\nonumber\\
&&\times {\cal P}(t)\Big(\bar{A}_0(t)+\bar{D}_{0}(t) {\cal P}(t) \Big).
\end{eqnarray}
By combining (\ref{z16}) with (\ref{sec2}) and taking similar transformations to $\Upsilon_2(t)$, we have
\begin{eqnarray}
0&=&\dot{{\cal P}}(t)+{\cal P}(t)A_0(t)+A_0'(t){\cal P}(t)+{\cal P}(t) D_0(t) {\cal P}(t)\nonumber\\
&&+\Big[\bar{A}_0(t)+\bar{D}_0(t) {\cal P}(t)
\Big]'\Big[I-{\cal P}(t) \bar{F}_0(t)\Big] {\cal P}(t) \nonumber\\
&&\Big[\bar{A}_0(t)
+\bar{D}_0(t) {\cal P}(t) \Big] -\Gamma_2(t)'\Upsilon_2(t) \Gamma_2(t), \label{TJ1}
\end{eqnarray}
with ${\cal P}(T)=P_1(T)+P_2(T)$ and $\Gamma_2(t)$ and $\Upsilon_2(t)$ are given as in the above. Note that Riccati equation (\ref{TJ1}) is exactly the same as  (\ref{p1}) with the arbitrary terminal values. Recall that $\Upsilon_1(t)=0$ and $\Gamma_1(t)\not=0$ for arbitrary terminal value $P_1(T)$, thus $\Upsilon_2(t)=0$ and $\Gamma_2(t)\not=0$ for arbitrary terminal value ${\cal P}(T)$.  This implies that where is no $u_2(t)$ such that (\ref{c59}) hold, thus Problem (IR-LQ) is unsolvable.
The proof is completed now. \hfill $\blacksquare$

Now we have the following Observation.

\begin{observation}
The optimization procedure for Problem (IR-LQ) is completed in finite layers, i.e., at most $m$ layer where $m$ is the dimension of control input $u(t)$. The optimization will stop  in one of the following two cases:
 \begin{itemize}
 \item $Range \Big(\Gamma_k(t)\Big)\subseteq Range \Big(\Upsilon_k(t)\Big)$;
 \item $\Gamma_k(t)\not=0$ and $\Upsilon_k(t)=0$.
\end{itemize}
\end{observation}

{\em Proof:}  By using Theorem \ref{JNB1} and noting that $dim(u_{k}(t))<dim(u_{k-1}(t))$, we know that  the proof follows directly.


\section{Open-loop and closed-loop solutions}

According to Theorem \ref{theorem5}, if $Range \Big(\Gamma_i(t)\Big)\subseteq Range \Big(\Upsilon_i(t)\Big),$ then Problem (IR-LQ) is solvable if and only if there exists $z_{i}(t)$ such that $\sum_{j=1}^iP_j(T)x(T)=0$ where $x(t)$ obeys the following dynamic:
\begin{eqnarray}
dx(t)&=&\Big[\Big(M_i(t)-H_{i}(t)\Upsilon_i^{\dag}(t)\Gamma_i(t)\Big)x(t)\nonumber\\
&&+H_{i}(t)[I-\Upsilon_i^{\dag}(t)\Upsilon_i(t)]z_{i}(t)\Big]dt+\Big[\Big(\bar{M}_i(t)\nonumber\\
&&-\bar{H}_{i}(t)\Upsilon_i^{\dag}(t)\Gamma_i(t)\Big)x(t)+\bar{H}_{i}(t)[I\nonumber\\
&&-\Upsilon_i^{\dag}(t) \Upsilon_i(t)]z_{i}(t)\Big]dw(t). \nonumber
\end{eqnarray}
\begin{remark}
In particularly, it should be pointed out that if $\Upsilon_i(t)$ is invertible then term  $z_{i}(t)$ will disappear and the dynamic system becomes
\begin{eqnarray}
dx(t)&=&\Big[\Big(M_i(t)-H_{i}(t)\Upsilon_i^{-1}(t)\Gamma_i(t)\Big)x(t)\Big]dt\nonumber\\
&&+\Big[\Big(\bar{M}_i(t)-\bar{H}_{i}(t)\Upsilon_i^{-1}(t)\Gamma_i(t)\Big)x(t)\Big]dw(t).\nonumber\\ \label{dynamic1}
\end{eqnarray}
Then Problem (IR-LQ) is solvable if and only if $\sum_{i=j}^iP_j(T)x(T)=0$ where $x(t)$ obeys the above dynamic (\ref{dynamic1}).
\end{remark}

Following the denotations in last subsection, the above equation can be simply denoted as
\begin{eqnarray}
dx(t)&=&\Big[A_i(t)x(t)+B_{i}(t)u_{i}(t)\Big]dt+\Big[\bar{A}_i(t)x(t)\nonumber\\
&&+\bar{B}_{i}(t)u_{i}(t)\Big]dw(t).\label{m17}
\end{eqnarray}
Recalling Theorem 1 of \cite{pengshige2}, it has been obtained that the system (\ref{m17}) is exactly controllable if and only if $\bar{B}_{i}(t)$ has full row rank.
In this case, there exits a matrix $\Lambda_i(t)$ such that
\begin{eqnarray}
\bar{B}_{i}(t)\Lambda_{i}(t)=\left[
                       \begin{array}{cc}
                         I & 0 \\
                       \end{array}
                     \right].\nonumber
\end{eqnarray}
Denote
\begin{eqnarray}
\Lambda_{i}^{-1}(t)u_{i}(t)&=&\left[
                    \begin{array}{c}
                      \mu_{i}(t) \\
                      \nu_{i}(t) \\
                    \end{array}
                  \right],\nonumber\\
\bar{y}_i(t)&=&\bar{A}_i(t)x(t)+\mu_{i}(t),\nonumber\\
B_{i}(t)\Lambda_i(t)&=&\left[
             \begin{array}{cc}
               G_{i1}(t) & G_{i2}(t) \\
             \end{array}
           \right].\nonumber
\end{eqnarray}
Rewriting (\ref{m17}) yields that
\begin{eqnarray}
dx(t)&=&-\Big(-[A_i(t)-G_{i1}(t)\bar{A}_i(t)]x(t)-G_{i1}(t)\bar{y}_i(t)\nonumber\\
&&-G_{i2}(t)\nu_{i}(t)\Big)dt+\bar{y}_i(t)dw(t)\label{c23}
\end{eqnarray}

\begin{theorem}\label{theorem2}
If $Range\Big(\Gamma_i(t)\Big)\subseteq Range \Big(\Upsilon_i(t)\Big)$ and $\bar{B}_i(t)$ has full-row rank,
then there exists an open-loop/closed-loop solution to Problem (  IR-LQ) if and only if there exists an
open-loop/closed-loop solution such that there exists an open-loop/closed-loop controller $\nu_{i}(t)$ such that $\sum_{j=1}^iP_j(T)x(T)=0$ where $x(t)$ is given by (\ref{c23}).
In particular, it follows that

1). The open-loop solution can be given by
$u_{i}(t)=\Lambda_i(t)\left[
                    \begin{array}{c}
                      \mu_{i}(t) \\
                      \nu_{i}(t) \\
                    \end{array}
                  \right]$
with \begin{eqnarray}
\nu_{i}(t)=-G_{i2}'(t)\Phi_{i}(t)G_{i}[t_0,T]^{-1}T_i'(t_0)x_0,\label{m13}
\end{eqnarray}
if the following Gramian matrix is invertible:
\begin{eqnarray}
G_{i}[t_0,T]&=&E\Big[\int_{t_0}^T\Phi_{i}'(s)G_{i2}(s)G_{i2}'(s)\Phi_{i}(s)ds\Big]
\end{eqnarray}
while $\Phi_{i}(t)$ satisfies
\begin{eqnarray}
d\Phi_{i}(t)&=&-[A_i(t)-G_{i1}(t)\bar{A}_i(t)]'\Phi_{i}(t)dt\nonumber\\
&&-G_{i1}'(t)\Phi_{i}(t)dw(t),\nonumber\\
\Phi_{i}(t_0)&=&I.\nonumber
\end{eqnarray}

2). The closed-loop solution can be given by
$z_{i}(t)=K_i(t)x(t)$ where $K(t)$ satisfying
\begin{eqnarray}
[A_i(t)-G_{i1}(t)\bar{A}_i(t)]+G_{i2}(t)K(t)=\frac{1}{t-T}.\nonumber
\end{eqnarray}
if it holds that
\begin{eqnarray}
Range\Big(G_{i2}(t)\Big)&\subseteq& Range\Big(\frac{1}{t-T}-[A_i(t)\nonumber\\
&&-G_{i1}(t)\bar{A}_i(t)]\Big)\nonumber
\label{c24}
\end{eqnarray}
\end{theorem}
\emph{Proof.}
1). From (\ref{c23}), it yields that
\begin{eqnarray}
x(t)&=&x(T)-E\Big[{\Phi_{i}'(t)}^{-1}\int_t^T\Phi_{i}'(s)G_{i1}'(s)\nu_{i1}(s)ds|\mathcal{F}_t\Big].\nonumber
\end{eqnarray}
This implies that
\begin{eqnarray}
x(t_0)&=&x(T)-E\Big[\int_{t_0}^T\Phi_{i}'(s)G_{i1}'(s)\nu_{i1}(s)ds\Big].\label{m18}
\end{eqnarray}
Let $\nu_{i1}$ satisfy (\ref{m13}), then (\ref{m18}) can be reformulated as
\begin{eqnarray}
x(t_0)&=&x(T)+x_0.
\end{eqnarray}
Accordingly, $x(T)=0.$ This implies that $\sum_{j=1}^iP_j(T)x(T)=0.$

2). If (\ref{c24}) holds, then there exists a matrix $K(t)$ such that
\begin{eqnarray}
[A_i(t)-G_{i1}(t)\bar{A}_i(t)]+G_{i2}(t)K(t)=\frac{1}{t-T}.\nonumber
\end{eqnarray}
Let $z_{i}(t)=K(t)x(t)$ and $\bar{y}_i(t)=0,$ then (\ref{c23}) becomes
\begin{eqnarray}
dx(t)&=&\frac{1}{t-T}x(t)dt.
\end{eqnarray}
By solving the above equation, we have
\begin{eqnarray}
x(T)=0.\nonumber
\end{eqnarray}
Thus, $\sum_{j=1}^iP_j(T)x(T)=0.$
The proof is completed. \hfill $\blacksquare$


\section{Examples}\label{ex}

Consider the system (\ref{r22}) and the corresponding cost function in the first example in Section I. The Riccati equation is as
\begin{eqnarray}
\dot{P}(t)=0,~P(1)=1,\nonumber
\end{eqnarray}
this implies that $P(t)=1,~t\in[0,1]$
Thus
\begin{eqnarray}
R+\bar{B}'P(t)\bar{B}=\left[
            \begin{array}{cc}
                          1 & -1 \\
                          -1 & 1 \\
                        \end{array}
                      \right]\nonumber\\
B'P(t)+\bar{B}'P(t)\bar{A}=\left[
                             \begin{array}{c}
                               1 \\
                               1 \\
                             \end{array}
                           \right]
\end{eqnarray}
It is obvious that the solution $P(t)$ is nonregular. The solvability is equivalent to the solvability of the following FBSDEs:
\begin{eqnarray}
dx(t)&=&\left[
          \begin{array}{cc}
            1 & 1 \\
          \end{array}
        \right]
z(t)dt-\bar{\Theta}(t)dw(t),\nonumber\\
d\Theta(t)&=&-\left[
          \begin{array}{cc}
            1 & 1 \\
          \end{array}
        \right]z(t)dt+\bar{\Theta}(t)dw(t),\nonumber\\
0&=&\left[
      \begin{array}{c}
        1 \\
        1 \\
      \end{array}
    \right]
x(t)+\left[
      \begin{array}{c}
        1 \\
        1 \\
      \end{array}
    \right]
\Theta(t)+\left[
      \begin{array}{c}
        1 \\
        1 \\
      \end{array}
    \right]\bar{\Theta}(t),\nonumber
\end{eqnarray}
with $\Theta(T)=0.$
Following Theorem \ref{the2}, let $P_1(t)=-P(t)=-1,$ if $P_1(T)x(T)=-x(T)=0,$ then the optimization problem is solvable.
In this case, together with Theorem \ref{theorem2}, the open-loop solution is given by
\begin{eqnarray}
u(t)
&=&-\left[
      \begin{array}{c}
        1 \\
        1 \\
      \end{array}
    \right]\frac{x_0}{2(1-t_0)}.\nonumber
\end{eqnarray}
The closed-loop solution is given by
\begin{eqnarray}
u(t)
&=&\left[
      \begin{array}{c}
        1 \\
        1 \\
      \end{array}
    \right]\frac{x(t)}{2(t-1)}.\nonumber
\end{eqnarray}


\section{Conclusions}

In this paper we have solved the irregular-LQ optimal control problem by proposing the  multiple-layer optimization approach, which is new to our best of knowledge. The controllers including the open-loop and closed-loop are designed based on Riccati equations and terminal constraint of state. The results show the essential difference between the irregular-LQ control and the regular-LQ control which lies in that the former must guarantee the terminal constraint of state ($Mx(T)=0$).


\end{document}